%% file: 2005-23.tex
\documentclass{gtart_h}  

\input gtoutput

\lognumber{489}

\received{22 August 2004}
\volumenumber{9}\papernumber{23}\volumeyear{2005}
\pagenumbers{991}{1042}   
\published{1 June 2005}
\accepted{13 April 2005}
\proposed{Peter Ozsvath}
\seconded{Walter Neumann, Tomasz Mrowka}

\usepackage{euscript}
\usepackage{amsmath}
\usepackage{amssymb}

\theoremstyle{definition}
\swapnumbers
\newtheorem{h}[subsection]{}
\newtheorem{hh}[subsubsection]{}

\newtheorem{exss}[subsubsection]{\,Example}
\newtheorem{defs}[subsection]{\,Definition}

\theoremstyle{plain}
\newtheorem{thms}[subsection]{\,Theorem}
\newtheorem{props}[subsection]{\,Proposition}
\newtheorem{cors}[subsection]{\,Corollary}
\newtheorem{lems}[subsection]{\,Lemma}
\newtheorem{thmss}[subsubsection]{\,Theorem}
\newtheorem{propss}[subsubsection]{\,Proposition}

\newtheorem{lemss}[subsubsection]{\,Lemma}

\newcommand{\Char}{\operatorname{Char}}
\newcommand{\Spin}{\operatorname{Spin}}
\newcommand{\spin}{\operatorname{spin}}

\DeclareMathOperator{\Hom}{{\rm Hom}}
\DeclareMathOperator{\rank}{{\rm rank}}
\newcommand{\ssw}{\textbf{sw}}
\newcommand{\et}{\EuScript{T}}
\newcommand{\bA}{{\mathbb A}}
\newcommand{\bB}{{\mathbb B}}
\newcommand{\bH}{{\mathbb H}}
\newcommand{\bms}{\mbox{\boldmath$s$}}
\newcommand{\xo}{o}
\def\C{\mathbb C}
\def\Q{\mathbb Q}
\def\R{\mathbb R}
\def\Z{\mathbb Z}
\def\N{\mathbb N}
\def\im{{\rm Im}}
\newcommand{\lp}{{l}}
\newcommand{\ev}{\varepsilon}
\newcommand{\cale}{{\mathcal E}}
\newcommand{\caly}{{\mathcal Y}}
\newcommand{\calx}{{\mathcal X}}
\newcommand{\calv}{{\mathcal V}}
\newcommand{\calt}{{\mathcal T}}
\newcommand{\calj}{{\mathcal J}}

\begin{document}

\title{On the  Ozsv\'ath--Szab\'o invariant of negative\\definite
plumbed 3--manifolds}
\covertitle{On the  Ozsv\noexpand\'ath--Szab\noexpand\'o invariant of negative definite
plumbed 3--manifolds}
\asciititle{On the  Ozsvath-Szabo invariant of negative definite
plumbed 3-manifolds}
\shorttitle{On the Ozsv\'ath--Szab\'o  invariant of plumbed 3--manifolds}

\author{Andr\'as N\'emethi}
\coverauthors{Andr\noexpand\'as N\noexpand\'emethi}
\asciiauthors{Andras Nemethi}
\address{Department of Mathematics\\Ohio State University\\Columbus, OH
43210, USA}
\secondaddress{R\'enyi Institute of Mathematics\\Budapest,
Hungary}
\asciiaddress{Department of Mathematics\\Ohio State University\\Columbus, OH
43210, USA\\and\\Renyi Institute of Mathematics\\Budapest,
Hungary}

\gtemail{\mailto{nemethi@math.ohio-state.edu}, \mailto{nemethi@renyi.hu}}
\asciiemail{nemethi@math.ohio-state.edu, nemethi@renyi.hu}
\urladdr{http://www.math.ohio-state.edu/~nemethi/}

\begin{abstract}
The main goal of the present article is the computation of the
Heegaard Floer homology introduced by Ozsv\'ath and Szab\'o
for a family of plumbed rational homology 3--spheres.
The main motivation is the study
of the Seiberg--Witten type invariants of links of normal surface singularities.
\end{abstract}

\asciiabstract{%
The main goal of the present article is the computation of the
Heegaard Floer homology introduced by Ozsvath and Szabo for a family
of plumbed rational homology 3-spheres.  The main motivation is the
study of the Seiberg-Witten type invariants of links of normal surface
singularities.}


\asciikeywords{3-manifolds, Ozsvath-Szabo Heegaard Floer homology,
Seiberg-Witten invariants, Seifert manifolds, Lens spaces,
Casson-Walker invariant, Q-homology spheres, Reidemeister--Turaev
torsion, normal surface singularities, rational singularities,
elliptic singularities} 

\maketitle

\section{Introduction}

The main goal of the present article is the computation of the
Ozsv\'ath--Szab\'o $\Z[U]$--module (Heegaard Floer homology)
$HF^+(M,[k])$ \cite{OSz,OSz7,OSzTr,OSzAB} for a family of plumbed
rational homology 3--spheres $M$ and any $\spin^c$--structure $[k]\in
\Spin^c(M)$. The author's main motivation is the study of
the Seiberg--Witten type invariants of links of normal surface
singularities.

More precisely, we consider that family of connected, negative
definite plumbing graphs -- we call them {\em AR}-graphs (see
Section~\ref{sec8}) -- which satisfy the following property: there exists a
vertex such that decreasing the decoration (the Euler number) of that
vertex we get a rational graph (in the sense of Artin).  This
class is surprisingly large: it contains the links of rational and
weakly elliptic singularities, the graphs considered by Ozsv\'ath
and Szab\'o in \cite{OSzP}, in particular, all the Seifert
manifolds (associated with  negative definite plumbing graphs).
For such plumbing graphs $\Gamma$ (and associated plumbed
3--manifolds $M(\Gamma)$) we prove (extending the main result of
\cite{OSzP}) the completely topological description
$HF^+(M(\Gamma),[k])=\bH^+(\Gamma,[k])$, where $\bH^+(\Gamma,[k])$
is the combinatorial $\Z[U]$--module introduced in \cite{OSzP}.

Moreover, for such graphs, we provide a precise combinatorial
formula (algorithm) for $\bH^+(\Gamma,[k])$. In order to do this, we
define a ``graded root'' $(R_k,\chi_k)$ associated with any
connected, negative definite plumbing graph $\Gamma$ and
characteristic element $k$.  This object connects in a mysterious
way two different types of properties, objects and invariants:
those coming from the Ozsv\'ath--Szab\'o (or Seiberg--Witten) theory
with those coming from algebraic geometry and singularity theory.
For example, its grading $\chi_k$ is in fact a Riemann--Roch formula,
which guides (modulo a shift) the absolute grading of
$HF^+(M,[k])$.  We believe  that $(R_k,\chi_k)$ is the right
object which guides the hierarchy of the topological types of
links of normal surface singularities (see, for example, Section
\ref{sec6});
and at the same time, it contains all the information about
$\bH^+(\Gamma,[k])$. In the body of the paper, for any {\em AR}
graph, we determine all these ``graded roots'', we read  from them
the combinatorial modules $\bH^+(\Gamma,[k])$, and we compute the
numerical invariants $\ssw^{OSz}_{M,[k]}$ which are the candidates
provided by the Ozsv\'ath--Szab\'o  theory for the Seiberg--Witten
invariants.

We exemplify the theory with a detailed discussion of the lens
spaces and Seifert manifolds. Additionally, in  these two  cases
we verify that $\ssw^{OSz}_{M,[k]}$ coincides  with the
Reidemeister--Turaev sign-refined torsion (function) normalized by
the Casson--Walker invariant -- a fact conjectured by Ozsv\'ath and
Szab\'o.

Examples show that the above results about $\bH^+(\Gamma,[k])$
cannot be extended for a substantially larger class of plumbing
graphs (than the {\em AR} graphs) -- unless one modifies the
definition of the combinatorial object $\bH^+$.

Although the paper is more topological/combinatorial, in some
places we emphasize the relevance of the results in the theory of
normal surface singularities, especially in the light of
\cite{SWI} and \cite{NL}.

The reader is invited to read Section~\ref{quid}, which contains  more
details about the guiding problems (and the main results) of the present
paper.

The author thanks P\'eter Ozsv\'ath and Zolt\'an Szab\'o for their help,
advice and encouragement.
The author is partially supported by NSF grant DMS-0304759
and OTKA grants 42769 and 46878.

\section{Preliminaries}

\begin{h}\textbf{Plumbing graphs}\qua\label{3.1}
We consider an oriented  plumbed 3--manifold $M=M(\Gamma)$ obtained
from the graph $\Gamma$ by plumbing construction (with all the
edge-signs equal to +; see \cite[Section~8]{NP}). We will
assume that $M$ is a {\em rational homology sphere}. This is
equivalent to $\Gamma$ being a disjoint union of trees,
all the genera of the corresponding Riemann surfaces involved in
the plumbing construction being zero, and the bilinear form $B$
associated with $\Gamma$ (see below) being non-degenerate.
Additionally, we will assume that $\Gamma$ is {\em connected}  and
$B$ {\em is negative definite} (for motivation see Example~\ref{3.2},
although most of the arguments work for non-connected $\Gamma$ as
well). In fact, the plumbing construction provides an oriented
4--manifold $\widetilde{X}$ whose oriented boundary $\partial
\widetilde{X}$ is exactly $M$. Let $-M$ denote $M$ with the opposite
orientation.

The vertices of $\Gamma$ will be
 indexed by $\calj$. Each vertex $j\in\calj$  has a decoration
$e_j$, which is the Euler number of the corresponding disc-bundle
(or $S^1$--bundle) used in the plumbing construction. Let $L$ be
the free $\Z$--module of rank  $s:=\#\calj$ with a fixed basis
$\{b_j\}_{j\in \calj}$.

In fact, $\Gamma$ is codified by the bilinear  form
$B=(b_i,b_j)_{i,j\in\calj}$, defined by $(b_j,b_j)=e_j$; and  for
$i\not=j$ one takes
 $(b_i,b_j)=1$ if the corresponding vertices $i,j$ are connected
by an edge in  $\Gamma$, otherwise $(b_i,b_j)=0$.
For any vertex $j\in\calj $ of $\Gamma$ we denote by $\delta_j$ the number
of adjacent vertices.

A cycle $x=\sum_j  n_jb_j\in L$ is  called \emph{effective}, denoted by $x\geq 0$, if $n_j\geq 0$
for any $j$.  We write $x\geq y$ if  $x-y\geq 0$. $x>0$ means $x\geq 0$ but $x\not=0$.
We define the support $|x|$ of $x$ by
$\bigcup_{j:n_j\not=0}\{b_j\}$.  If $x_i=\sum_jn_{j,i}b_j$ for $i=1,2$, then
$\min\{x_1,x_2\}:=\sum_j\min\{n_{j,1},n_{j,2}\}b_j$.
We write $x^2$ to mean  $(x,x)$.

Set $L'=\Hom_\Z(L,\Z)$, and consider the
exact sequence $0\to L\stackrel{i}{\to} L'\to H\to 0$, where
$i(x)=(x,\cdot )$. Then one has the identifications
$L\approx H_2(\widetilde{X},\Z),
\ L'\approx H_2(\widetilde{X},M,\Z)$ and $H\approx H_1(M,\Z)$.

Sometimes it is convenient to identify  the lattice $L'$ with a
sub-lattice of $L_\Q=L\otimes \Q$:   any $\alpha\in L'$ is identified with
$y_\alpha\in L_\Q$  which satisfies  $\alpha(x)=(y_\alpha , x)$ for any $x\in L$.
Here  $(\cdot,\cdot)$  is  the natural extension of the
form $B$  to $L_\Q$ (denoted in the same way).

In Sections \ref{sec10} and \ref{sec11} we will also use the dual base $\{g_j\}_{j\in\calj}$ of $L'$
(ie $(g_j,b_i)=\delta_{ij}$ for any $i,j\in\calj$).
If $L'$ is identified with a sub-lattice of $L_\Q$
then $g_j$ is the $j^{th}$ column of $B^{-1}$.
Since $B$ is negative definite, $ -g_j\geq 0$ for any $j\in \calj$.

Clearly, the order $|H|$ of the group $H=H_1(M,\Z)$  equals
$|\det(B)|$.
\end{h}

\begin{exss}(See \cite{NP})\qua\label{3.2} 
Let $M$ be  the link of a normal surface
singularity $(X,0)$.  Assume that $\widetilde{X}\to X$ is a good resolution of $(X,0)$,
ie  $X$ is a sufficiently small Stein representative of $(X,0)$ and
$\widetilde{X}\to X$ is a resolution with normal crossing  exceptional divisor.
Let  $\Gamma$ denote the dual resolution graph. Then  $\partial \widetilde{X}=M$,  and
 $(\widetilde{X},\partial \widetilde{X})$  can be identified with the output of the plumbing
construction provided by  $\Gamma$. (In this case, the elements
$b_j$ of $L=H_2(\widetilde{X},\Z)$ correspond to the classes of the
irreducible exceptional divisors.)

Recall, that a plumbed 3--manifold $M$ can be realized as the link
of some  normal surface singularity if and only if its graph
$\Gamma$ is {\em connected}  and the corresponding form $B$ {\em
is negative definite} \cite{GRa}.\end{exss}

\begin{h}\textbf{$\Z[U]$--modules}\qua
\label{zu} 
Later we will use the following notation.
Consider the graded $\Z[U]$--module $\Z[U,U^{-1}]$, and (following
\cite{OSzP}) denote by
 $\calt_0^+$ its quotient by the submodule  $U\cdot \Z[U]$.
This has a grading such that $\deg(U^{-d})=2d$ (for $d\geq
0$). Similarly, for any $n\geq 1$, define the  graded module
$\calt_0(n)$ by the quotient of $\Z\langle U^{-(n-1)},
U^{-(n-2)},\ldots, 1,U,\ldots \rangle$ by $U\cdot \Z[U]$ (with the
same grading). Hence, $\calt_0(n)$, as a $\Z$--module, is the free
$\Z$--module $\Z \langle 1,U^{-1},\ldots,U^{-(n-1)} \rangle$
(generated by $1,U^{-1},\ldots,U^{-(n-1)}$), and has finite
$\Z$--rank $n$.

More generally, for any graded $\Z[U]$--module $P$ with
$d$--homogeneous elements $P_d$, and  for any  $r\in\Q$,   we
denote by $P[r]$ the same module graded in such a way that
$P[r]_{d+r}=P_{d}$. Then set $\calt^+_r:=\calt^+_0[r]$ and
$\calt_r(n):=\calt_0(n)[r]$.
\end{h}

\begin{h}\textbf{Characteristic elements and $\Spin^c$--structures}\qua
\label{3.3}
Fix a plumbing graph $\Gamma$ as above. The set of characteristic
elements is defined by
$$\Char=\Char(\Gamma):=\{k\in L': \, k(x)+(x,x)\in 2\Z \ \mbox{for every $x\in L$}\}.$$
Clearly, for any fixed $k_0\in \Char$, $\Char=k_0+2L'$. There is a natural
action of $L$ on $\Char$ by $x*k:=k+2i(x)$ whose orbits are of type
$k+2i(L)$. Obviously, $H$ acts freely and transitively on the set of orbits
by $[l']*(k+2i(L)):=k+2l'+2i(L)$ (in particular,
they have the same cardinality).

If $\widetilde{X}$ is as in Section~\ref{3.1}, then  the first Chern class (of
the determinant line bundle associated with a given
$\spin^c$--structure) realizes an identification between the
$\spin^c$--structures $\Spin^c(\widetilde{X})$ on $\widetilde{X}$ and
$\Char\subset L'\approx H^2(\widetilde{X},\Z)$ (see
\cite[Section~2.4.16]{GS}). The restrictions to $M$ defines an
identification of the
$\spin^c$--structures $\Spin^c(M)$ of $M$ with the set of orbits of
$\Char$ modulo $2L$; and this identification is compatible with the
action of $H$ on both sets. In the sequel, we think about
$\Spin^c(M)$ by this identification, hence any $\spin^c$--structure
 of $M$ will be represented by an orbit  $[k]:=k+2i(L)\subset \Char$.

The $\spin^c$--structures $\Spin^c(M)$ and $\Spin^c(-M)$ are
 canonically identified.
\end{h}

\begin{h}\textbf{Invariants of the 3--manifold $M$}\qua
\label{invM}
In this note we will focus on the following set of  invariants.
\end{h}

\begin{hh} According to Turaev \cite{Tu5}, a choice of a
$\spin^c$--structure on $M$ is equivalent to a choice of an Euler
structure. For every $\spin^c$--structure $[k]$, he defines the
\emph{sign-refined Reidemeister--Turaev torsion}
$$\et_{M,[k]}=\sum_{h\in H}\et_{M,[k]}(h) h\in \Q[H]$$
determined by the Euler structure associated with $[k]$
(see \cite{Tu5}). It is convenient to  think about $\et_{M,[k]}$ as a
function $H\to\Q$ given by $h\mapsto \et_{M,[k]}(h)$. We will be
basically interested in $\et_{M,[k]}(1)$ where 1 denotes the
neutral element of the group $H$ (with the multiplicative
notation).

For plumbed 3--manifold $M$ (as in Section~\ref{3.1}), \cite[Section~5.7]{SWI} 
provides a combinatorial formula for $\et_{M,[k]}(1) $ in terms of
$\Gamma$ (involving regularized Fourier--Dedekind  sums).
\end{hh}

\begin{hh} We will use the notation $\lambda(M)$ for
the \emph{Casson--Walker invariant of $M$} (normalized as in
\cite[Section~4.7]{Lescop}). For plumbed manifolds $M$, it  has a
combinatorial formula from $\Gamma$ (provided by A. Ra\c{t}iu in
his dissertation; and it can be deduced from the surgery formulas
of \cite{Lescop} as well; see also  \cite[Section~5.3]{SWI}):
$$-\frac{24}{|H|}\,\lambda(M)=\sum_je_j+3s+\sum_j(2-\delta_j)\,
(B^{-1})_{jj}$$\end{hh}

\begin{hh}\label{knegyzet}  If $M$ is a plumbed
3--manifold as in Section~\ref{3.1}, it has a \emph{canonical characteristic
element} $K$ defined by $K(b_j)=-e_j-2$ for any $j$ (see also
\ref{4.2} for more comments about these ``adjunction formulas'').
The rational number $K^2+s$ does not depend on the (negative
definite) plumbing graph, it is an invariant of $M$. It can be
computed from $\Gamma$ as follows (see \cite[Section~5.2]{SWI}):
$$K^2+s=\sum_je_j+3s+2 +\sum_{i,j}(2-\delta_i)
(2-\delta_j)(B^{-1})_{ij}.$$\end{hh}

\begin{hh}\textbf{The Heegaard Floer homology $HF^+(M,[k])$}\qua For any
oriented rational homology 3--sphere   $M$ and $[k]\in \Spin^c(M)$,
the {\em Ozsv\'ath--Szab\'o $\Z[U]$--module} (or, the \emph{absolutely
graded Heegaard Floer homology}) $HF^+(M,[k])$ was introduced in
\cite{OSz} (see also the long list of recent articles of
Ozsv\'ath and Szab\'o). This has a $\Q$--grading compatible with
the $\Z[U]$--action, where $\deg(U)=-2$.

Additionally, $HF^+(M,[k])$ has an (absolute) $\Z_2$--grading;
$HF^+_\text{even}(M,[k])$, respectively
$HF^+_\text{odd}(M,[k])$ denote the part of
$HF^+(M,[k])$ with the corresponding parity.

By the general theory, for any $\spin^c$--structure $[k]$ and $M$ as above,
one has a graded $\Z[U]$--module isomorphism
$$HF^+(M,[k])=\calt^+_{d(M,[k])}\oplus HF^+_\text{red}(M,[k]),$$
where $HF^+_\text{red}$ has a finite $\Z$--rank and an induced
(absolute) $\Z_2$--grading (and $d(M,[k])$ can also be defined via this
isomorphism).

From the above data one can extract two crucial numerical invariants:
$d(M,[k])$ and
$$\chi(HF^+(M,[k])):=\rank_\Z HF^+_\text{red,even}(M,[k])
-\rank_\Z HF^+_\text{red,odd}(M,[k]).$$
With respect to the change of orientation they  behave as follows:
$$d(M,[k])=-d(-M,[k]) \ \  \mbox{and} \ \
\chi(HF^+(M,[k]))=-\chi(HF^+(-M,[k])).$$ Notice that one can
recover $HF^+(M,[k])$ from $HF^+(-M,[k])$ by  a standard procedure
described by Ozsv\'ath and Szab\'o.
\end{hh}

\begin{hh}\label{bhplus} \textbf{The combinatorial module
$\bH^+(\Gamma)$}\qua In the case of a plumbed 3--manifold $M$
associated with the graph $\Gamma$, one would like to obtain a
completely combinatorial description of $HF^+(-M,[k])$ from
$\Gamma$. In \cite{OSzP}, Ozsv\'ath and Szab\'o introduced the
following combinatorial graded $\Z[U]$--module as a candidate.

Let $\bH^+(\Gamma)$ be the set of those finitely supported functions
$\phi\co \Char\to \calt_0^+$ which satisfy the following property.
For any characteristic element $k$ and base vector $b_j$
write $2n=k(b_j)+(b_j,b_j)$. Then if $n\geq 0$ then one  requires
$$U^n\cdot \phi(k+2i(b_j))=\phi(k);$$
while if $n\leq 0$, then
 $$ \phi(k+2i(b_j))=U^{-n}\cdot \phi(k).$$

For any $\spin^c$--structure  $[k]=k+2i(L)$, let $\bH^+(\Gamma,[k])$
be the subset of those maps $\phi$ which are supported on $k+2i(L)$.
Then one has a direct sum splitting $\bH^+(\Gamma)=\bigoplus_{[k]}\,
\bH^+(\Gamma,[k])$.

On $\bH^+(\Gamma)$ one defines the following (rational)
grading. One says that a map
$\phi\in \bH^+(\Gamma)$ is homogeneous of degree $d\in \Q$ if for each
characteristic  vector $k$  with $\phi(k)\not=0$, $\phi(k)\in \calt_0^+$
is a homogeneous element with
$$\deg(\phi(k))-\frac{k^2+s}{4}=d.$$
 Ozsv\'ath and Szab\'o in \cite{OSzP} proved the following theorem.
 \end{hh}

\begin{thmss}
\label{OSZPI}
Assume that
$\Gamma$  has at most one vertex $j\in \calj$ where the inequality
$-e_j\geq \delta_j$ fails. Then for any $\spin^c$--structure $[k]$,
$HF^+_\text{odd}(-M,[k])=0$ and $HF^+_\text{even}(-M,[k])=
\bH^+(\Gamma,[k])$.
\end{thmss}

\begin{h}\textbf{The main guiding problems of the article}\qua
\label{quid}
We summarize in short those conjectures and  problems   which have
determined our interest in this subject. This is an interesting
mixture of the topology of 3--manifolds and theory of normal
surface singularities.
\end{h}

\begin{hh}\textbf{Problem I}\qua
First of all, as we already mentioned, for any plumbed manifold $M$ one
wishes to determine $HF^+(M,[k])$ in a combinatorial way from $\Gamma$.

Our goal (realized in Section~\ref{isohh}) is to extend the isomorphism
\ref{OSZPI} for
 a larger class of plumbing graphs, namely for the {\em AR} graphs
(introduced in section~\ref{sec8}). This class includes all the links of
rational and elliptic singularities, contains the class considered
by Ozsv\'ath and Szab\'o in Theorem~\ref{OSZPI}, in particular, all the
Seifert manifolds with negative orbifold Euler number. (Examples
show that the statement of Theorem~\ref{OSZPI} cannot be extended for a
substantially larger class of graphs; see Section~\ref{push}.)
\end{hh}

\begin{hh}\textbf{Problem II}\qua
Even though  $\bH^+(\Gamma,[k])$ is
combinatorial, its computation is not trivial at all. In the body
of the paper we provide precise algorithm for this computation
(valid for any {\em AR} graph).

The main idea comes from the technique of computation sequences
used in singularity theory (see \cite{Laufer72} and
\cite{Laufer77}). In fact, as a starting invariant, we will
construct a {\em graded root} $R_k$ for any plumbing graph
$\Gamma$ and for any characteristic element $k$. Surprisingly,
this object is able to guide (and connect) two very different
theories. First, from the point of view of singularity theory, it
seems that it is the right object which controls the hierarchy and
classification of surface singularities. For example, it characterizes
nicely the rational and elliptic singularities (see Section~\ref{sec6}).
Also,  it provides optimal topological upper bounds for some
analytic invariants (see Problem IV below and Section~\ref{pr7b}). On the
other hand, from $R_k$ one can read easily $\bH^+(\Gamma,[k])$.

Our discussion runs on three levels: we determine for any {\em AR}
graph the graded root $R_k$, which determines automatically
the combinatorial Ozsv\'ath--Szab\'o graded $\Z[U]$--modules
 $\bH^+(\Gamma,[k])$,
and finally we focus on the numerical invariants $d(M,[k])$ and
$\chi(HF^+(M,[k])$ as well.

One of the surprising byproducts of the discussion is Theorem~\ref{4.5},
which shows that the links of rational singularities correspond
exactly to the $L$--spaces in the sense of Ozsv\'ath--Szab\'o.
\end{hh}

\begin{hh}\textbf{Problem III}\qua
One of the goals of  Ozsv\'ath--Szab\'o homology
is to substitute the
(modified) Seiberg--Witten topological invariant of $(M,[k])$.
More precisely, it provides the numerical invariant
$$\ssw^{OSz}_{M,[k]}:=\chi(HF^+(M,[k]))-\frac{d(M,[k])}{2}$$
as a candidate for the Seiberg--Witten invariant.

On the other hand, the sign-refined Reidemeister--Turaev torsion
together with the Casson--Walker invariant also provides a
candidate for the Seiberg--Witten invariant (see also \cite{Nico5}):
$$\ssw^{TCW}_{M,[k]}:=-\et_{M,[k]}(1)+\frac{\lambda(M)}{|H|},$$
where $|H|$ is the order of $H$.

Our goal is to investigate the identity of these two candidates.
Notice that in the presence of $HF^+_\text{red,odd}(-M,[k])=0$,
this conjectured identity reads as
$$\et_{M,[k]}(1)-\frac{\lambda(M)}{|H|}=\rank_\Z \, HF^+_\text{red}(-M,[k])
-\frac{d(-M,[k])}{2}.$$
Recall that if $M$ is an integral homology sphere (ie $H=0$)
then this identity holds by \cite[Section~1.3]{OSzAB}.

Since for {\em AR} graphs the right hand side of the identity can be
determined from $\bH^+(\Gamma,[k])$, the above identity becomes a purely
combinatorial property of $\Gamma$. In the body of the paper we verify it
for any lens space and Seifert manifold (and any $\spin^c$--structure).
\end{hh}

\begin{hh}\textbf{Problem IV}\qua
In \cite{SWI},
L. Nicolaescu and the author stated the following conjecture.
Assume that $(X,0)$ is a normal surface singularity whose link $M$
is a rational homology sphere. Let $p_g$ be the geometric genus of
$(X,0)$ (ie $\dim_\C H^1(\widetilde{X},{\mathcal O}_{\widetilde{X}}) $,
where $\widetilde{X}$ is a resolution of $X$). Then the analytic
invariant $p_g$ has the following ``optimal'' topological upper
bound:
$$p_g\leq \et_{M,[K]}(1)-\frac{\lambda(M)}{|H|}-\frac{K^2+s}{8}$$
Moreover, if $(X,0)$ is $\Q$--Gorenstein, then above one has equality.

Notice that  in the presence of $HF^+_\text{red,odd}(-M,[K])=0$, and
using the principle in Problem III above, this inequality can be
transformed into another conjectured inequality
$$p_g\leq \rank_\Z (HF^+_\text{red}(-M,[K]))+ \frac{d(M,[K])}{2}-\frac{K^2+s}{8}$$
which conjecturally  becomes  equality for $\Q$--Gorenstein
singularities.

The point is that  the computation algorithm  of
$\bH^+(\Gamma,[K])$ of the present article will automatically
provide this inequality for any singularity with {\em
AR}-resolution graph. Moreover, we also prove the identity for
rational and elliptic Gorenstein singularities and all
singularities which admit a good $\C^*$--action (all of these are
$\Q$--Gorenstein).

In fact, in Corollary~\ref{pr7b}, we extend the above inequality for {\em
any} $\spin^c$ structure and any {\em AR} graph. In this way, the
(candidates for the)  Seiberg--Witten invariants (provided by the
Heegaard Floer theory)
 provide a (conjecturally optimal) topological upper bound for
$h^1({\mathcal L})$, where ${\mathcal L}$ is a holomorphic line
bundle on a resolution of $(X,0)$.

For the extension of the conjecture \cite{SWI} to arbitrary $\spin^c$
structures, see \cite{NL}.
For different definitions regarding
surface singularities, the reader is invited to consult \cite{Five}.

After the first version of the present article was completed, the
author realized that the above inequality -- and equality for
Gorenstein  singularities -- is not true for non-{\em AR}
singularities; for details see \cite{LMN}.
\end{hh}

\section{Graded roots}

\begin{h}\textbf{Preliminary remarks}\qua
\label{pr} We recall that the
Ozsv\'ath--Szab\'o $\Z[U]$--module $\bH^+(\Gamma,[k])$ is defined in
a combinatorial way from the plumbing graph $\Gamma$. Our goal is
to define an intermediate object, a graded root $R_k$ associated
with $\Gamma$ (and a characteristic element $k$). This will
contain all the needed information to determine the homological
object $\bH^+$, but it preserves also some additional, more subtle
topological information about $\Gamma$ (or, about $M$).

In this section we give the definition and first properties of the
(abstract) graded roots. The next section contains the construction
of the graded roots $R_k$ from the plumbing graphs $\Gamma$.
(Although both $\Gamma$ (the plumbing graph) and the constructed graded root
$R_k$ are ``connected trees'', they serve rather different roles.
For example, the edges of $\Gamma$ codify the corresponding gluings in the
plumbing, while  the edges of $R_k$ codify the $\Z[U]$--action.
We hope the terminology will not create any confusion.)
\end{h}

\begin{h}\textbf{Definitions}\qua
\label{2.1}
(1)\qua  Let $R$ be an infinite tree with
vertices $\calv$ and edges $\cale$. We denote by $[u,v]$ the edge with
 end-points $u$ and $v$.  We say that $R$ is a {\em graded root}
with grading $\chi\co\calv\to \Z$ if
\begin{enumerate}
\item[(a)] $\chi(u)-\chi(v)=\pm 1$ for any $[u,v]\in \cale$;

\item[(b)] $\chi(u)>\min\{\chi(v),\chi(w)\}$ for any $[u,v],\
[u,w]\in\cale$, $v\neq w$;

\item[(c)] $\chi$ is bounded below, $\chi^{-1}(k)$ is finite for any $k\in\Z$,
and $\#\chi^{-1}(k)=1$ if $k$ is sufficiently large.
\end{enumerate}

(2)\qua We say that $v\in\calv $ is a {\em local minimum
point} of the graded root $(R,\chi)$ if $\chi(v)<\chi(w)$ for any edge
$[v,w]$.

(3)\qua If $(R,\chi)$ is a graded root, and $r\in \Z$, then
we denote by $(R,\chi)[r]$ the same $R$ with the new grading
$\chi[r](v):=\chi(v)+r$. (This can be generalized for any $r\in
\Q$ as well.)
\end{h}

\begin{h}\textbf{Notation and remarks}\qua
\label{2.2}
(1)\qua For a vertex $v$ set $\delta_v:=\#\{[v,w]\in \cale\}$.
One can verify that the set of vertices  $\calv_1:=\{v\in \calv:\delta_v=1\}$
are exactly the local minimum points of $\chi$,  and $\#\calv_1<\infty$.

(2)\qua A geodesic path connecting two vertices is {\em monotone } if
$\chi$ restricted  to the set of vertices on the path is strict monotone.
If a vertex $v$ can be connected by another vertex $w$ by a monotone geodesic
and $\chi(v)>\chi(w)$, then we say that $v$ dominates $w$, and we write
$v\succ w$. $\succ$ is an ordering of $\calv$. For any pair $v,w\in\calv$
there is a unique $\succ$--minimal vertex $\sup(v,w)$ which dominates both.
\end{h}

\begin{h}\textbf{Examples}\qua
\label{2.3}
(1)\qua For any integer $n\in\Z$, let  $R_n$ be the tree with $\calv=\{v^{k}
\}_{ k\geq n}$ and $\cale=\{[v^{k},v^{k+1}]\}_{k\geq n}$. The grading
is $\chi(v^{k})=k$.

(2)\qua Let $I$ be a finite index set. For each $i\in I$ fix  an integer
$n_i\in \Z$; and for each pair $i,j\in I$ fix $n_{ij}=n_{ji}\in\Z$ with the
properties
\begin{enumerate}
\item[(i)] $n_{ii}=n_i$;
\item[(ii)] $n_{ij}\geq \max\{n_i,n_j\}$; and
\item[(iii)] $n_{jk}\leq \max\{n_{ij},n_{ik}\}$
\end{enumerate}
for any $ i,j,k\in I$.

For any $i\in I$ consider $R_i:=R_{n_i}$ with vertices $\{v_i^{k}\}$ and
edges $\{[v_i^{k},v_i^{k+1}]\}$, $(k\geq n_i)$.
In the disjoint union $\coprod_iR_i$,  for any pair $(i,j)$,
 identify $v_i^{k}$ and $v_j^{k}$,
respectively $[v_i^{k},v_i^{k+1}]$  and $[v_j^{k},v_j^{k+1}]$, whenever
$k\geq n_{ij}$. Write $\overline{v}_i^{k}$ for the class of $v_i^k$.
Then $\coprod_iR_i/_\sim$ is a graded root with
$\chi(\overline{v}_i^{k})=k$. It will be denoted by $R=R(\{n_i\},\{n_{ij}\})$.

Clearly $\calv_1(R)$
is a subset of $\{\overline{v}_i^{n_i}\}_{i\in I}$, and this last set
can be identified with $I$. $\calv_1(R)=I$ if  in (ii) all the inequalities
are strict. Otherwise all the indices $I\setminus \calv_1(R)$ are
superfluous, ie the corresponding $R_i$'s produce no additional vertices.

In fact, any graded root $(R',\chi')$ is isomorphic (in a natural sense)
with some $R(\{n_i\},\{n_{ij}\})$.
Indeed, set $I:=\calv_1(R')$,
$n_v:=\chi'(v)$ and $n_{uv}:=\chi'(\sup(u,v))$ for $u,v\in I$.

(3)\qua Any map $\tau\co\{0, 1,\ldots,l\}\to \Z$ produces starting data for
construction (2). Indeed, set $I=\{0,\ldots,l\}$, $n_i:=\tau(i)$
for $i\in I$, and $n_{ij}:=\max\{n_k\,:\, i\leq k\leq j\}$ for $i\leq j$.
Then  $\coprod_iR_i/_\sim $ constructed
in (2) using this data will be denoted by $(R_\tau,\chi_\tau)$.

For example, for $l=4$, take for the values of $\tau $:\  $-3,-1,-2,0$ and
$-2$ (respectively $-3,0,-2,-1$ and $-2$). Then the two graded roots are:

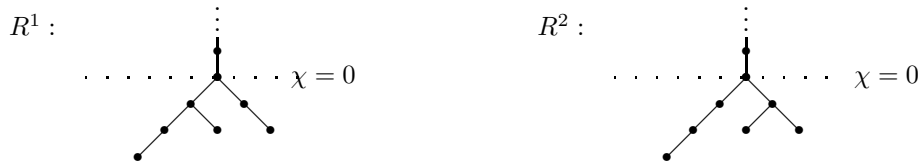
\begin{picture}(375,70)(15,0)\small
\linethickness{.5pt}
\put(100,40){\circle*{3}}
\put(100,50){\circle*{3}}
\put(90,30){\circle*{3}}
\put(80,20){\circle*{3}}
\put(70,10){\circle*{3}}
\put(100,40){\circle*{3}}
\put(100,20){\circle*{3}}
\put(110,30){\circle*{3}}
\put(120,20){\circle*{3}}
\put(100,40){\line(-1,-1){30}}
\put(100,40){\line(1,-1){20}}
\put(90,30){\line(1,-1){10}}
\put(100,40){\line(0,1){15}}
\put(100,65){\makebox(0,0){$\vdots$}}
\put(140,40){\makebox(0,0)[0]{$\chi=0$}}
\put(30,60){\makebox(0,0){$R^1:$}}
\qbezier[10](50,40)(90,40)(130,40)

\put(300,40){\circle*{3}}
\put(300,50){\circle*{3}}
\put(290,30){\circle*{3}}
\put(280,20){\circle*{3}}
\put(270,10){\circle*{3}}
\put(300,40){\circle*{3}}
\put(300,20){\circle*{3}}
\put(310,30){\circle*{3}}
\put(320,20){\circle*{3}}
\put(300,40){\line(-1,-1){30}}
\put(300,40){\line(1,-1){20}}
\put(310,30){\line(-1,-1){10}}
\put(300,40){\line(0,1){15}}
\put(300,65){\makebox(0,0){$\vdots$}}
\put(340,40){\makebox(0,0)[l]{$\chi=0$}}
\qbezier[10](250,40)(290,40)(330,40)
\put(230,60){\makebox(0,0){$R^2:$}}
\end{picture}

In fact, we can even start with the infinite index set $I=\N$, and
a function $\tau\co\N\to\Z$, provided that we assure that starting
from a bound $l$, all the contributions $R_i$ with $i>l$ are
superfluous. This happens, for example, if, for some $l$, $\tau(i+1)\geq
\tau(i)$ for  any $i\geq l$.
\end{h}

\begin{h}\textbf{Definition: the associated $\Z[U]$--modules}\qua
\label{2.6}
For a graded root $(R,\chi)$, let $\bH(R,\chi)$ (shortly $\bH(R)$)
be the set of functions $\phi\co\calv\to \calt^+_0$ with the
following property: whenever $[v,w]\in \cale$ with
$\chi(v)<\chi(w)$, then
\begin{equation*}
U\cdot \phi(v)=\phi(w).
\end{equation*}
Or, equivalently, for any $w\succ v$ one requires
\begin{equation*}
U^{\chi(w)-\chi(v)}\cdot \phi(v)=\phi(w).
\tag{$*$}\end{equation*}
Clearly $\bH(R)$ is a $\Z[U]$--module via $(U\phi)(v)=U\cdot \phi(v)$.
Moreover, $\bH(R)$ has a grading:
$\phi\in \bH(R)$ is homogeneous of degree $d\in\Z$ if for each $v\in\calv$
with $\phi(v)\not=0$, $\phi(v)\in\calt^+_0$ is homogeneous  of degree
$d-2\chi(v)$. Notice that in $(*)$ one has $2\chi(v)+\deg\phi(v)=
2\chi(w)+\deg\phi(w)$, hence $d$ is well-defined.

Notice that any $\phi$ as above is automatically finitely
supported.
\end{h}

\begin{hh} From the definitions, it is clear that
$\bH((R,\chi)[r])=\bH(R,\chi)[2r]$ for any $r\in\Z$.\end{hh}

\begin{propss}\label{2.7}
Let $(R,\chi)$ be a
graded root. Set $I:=\{v\in\calv: \delta_v=1\}$, and we order the
set $I$ as follows. The first element $v_1$ is an arbitrary vertex
with $\chi(v_1)=\min_{v}\chi(v)$. If $v_1,\ldots, v_k$ is already
determined, and $J:=\{v_1,\ldots,v_k\}\varsubsetneq I$, then let
$v_{k+1}$   be an arbitrary vertex in $I\setminus J$ with
$\chi(v_{k+1})=\min _{v\in I\setminus J}\chi(v)$. Let $w_{k+1}\in
\calv$ be the unique $\succ$--minimal vertex of $R$ which dominates
both $v_{k+1}$, and at least one vertex from $J$.  Then one has
the following isomorphism of $\Z[U]$--modules:
$$\bH(R,\chi)=\calt^+_{2\chi(v_1)}\oplus \bigoplus_{k\geq 2}\calt_{2\chi(v_k)}\big(
\chi(w_k)-\chi(v_k)\big)$$
In particular, with the notation  $m:=\min_v\chi(v)$ and
$$\bH_\text{red}(R,\chi):=\bigoplus_{k\geq 2}\calt_{2\chi(v_k)}\big(
\chi(w_k)-\chi(v_k)\big),$$ one has a canonical direct sum
decomposition
$$\bH(R,\chi)=\calt^+_{2m}\oplus \bH_\text{red}(R,\chi)$$
of graded $\Z[U]$--modules.
The
$\Z[U]$--module $\bH_\text{red}(R)$ has finite $\Z$--rank, with
$\bH_\text{red}(R)=0$ if and only if $\#I=1$ and $R=R_{\min\chi}$.
\end{propss}

\begin{proof}
The proof is elementary, and is left to the reader.
\end{proof}

\begin{h}\textbf{Examples}\qua
\label{2.9}
(a)\qua $\bH(R_n)=\calt_{2n}$.

(b)\qua The graded roots  $R^1$ and $R^2$ constructed in Example~\ref{2.3}(3)
are not isomorphic, but their $\Z[U]$--modules are:
$\bH(R^1)=\bH(R^2)= \calt^+_{-6}\oplus \calt_{-4}(1)\oplus
\calt_{-4}(2)$. Hence, in general, a  graded root carries more
information than  its $\Z[U]$--module.
\end{h}

\begin{cors}\label{2.10}
Let $(R_\tau,\chi_\tau)$ be a
graded root associated with  some  function $\tau\co\N\to \Z$ (see
Example~\ref{2.3}(3)) which satisfies $\tau(1)>\tau(0)$.  Then the
$\Z$--rank of $\bH_\text{red}( R_\tau,\chi_\tau)$ is
$$\rank_\Z \bH_\text{red}(R_\tau)=-\tau(0)+\min_{i\geq 0}\tau(i)+\sum_{i\geq 0}\,
\max\{ \tau(i)-\tau(i+1),0\}.$$ The summand $\calt^+_{2m}$ of
$\bH(T_\tau,\chi_\tau)$ has index $m=\min_{i\geq
0}\tau(i)=\min_v\chi_\tau(v)$.
\end{cors}

\begin{proof}
Use induction over $l$ (where $\tau\co\{0,\ldots,l\}\to\Z$).
\end{proof}

\section{Graded roots associated with  plumbing graphs}

\begin{h}
Fix a connected plumbing graph $\Gamma$ whose bilinear form
is negative definite (see Section~\ref{3.1}). In this section we will
construct a graded root $(R_k,\chi_k)$ associated with any
characteristic element $k$.
\end{h}

\begin{h}\textbf{The construction of $(R_k,\chi_k)$}\qua
\label{3.4}
Fix $k\in \Char$ and define $\chi_k\co L\to\Z$ by
$$\chi_k(x):=-(k(x)+(x,x))/2.$$
For any $n\in \Z$, we define a finite 1--dimensional simplicial complex
$\overline{L}_{k,\leq n}$ as follows. Its 0--skeleton is
$L_{k,\leq n}:=\{x\in L:\, \chi_k(x)\leq n\}$.
For each $x$ and $j\in \calj$ with $x,x+b_j\in L_{k,\leq n}$,
we consider a unique 1--simplex  with endpoints at $x$ and $x+b_j$
(eg, the segment $[x,x+b_j]$ in $L\otimes \R$).
We denote the set of connected components of $\overline{L}_{k,\leq n}$ by
$\pi_0(\overline{L}_{k,\leq n})$. For any $v\in \pi_0(\overline{L}_{k,\leq n})$, let
$C_v$ be the corresponding connected component of $\overline{L}_{k,\leq n}$.

Next, we define the  graded root  $(R_k,\chi_k)$
as follows. The vertices $\calv(R_k)$ are $\bigcup_{n\in \Z}
\pi_0(\overline{L}_{k,\leq n})$. The grading $\calv(R_k)\to\Z$, which we still denote by $\chi_k$, is
$\chi_k|\pi_0(\overline{L}_{k,\leq n})=n$.

If $v_n\in \pi_0(\overline{L}_{k,\leq n})$, and $v_{n+1}\in
\pi_0(\overline{L}_{k,\leq n+1})$, and $C_{v_n}\subset C_{v_{n+1}}$,
then $[v_n,v_{n+1}]$ is an edge of $R_k$. All the edges
$\cale(R_k)$ of $R_k$ are obtained in this way.
\end{h}

\begin{props}\label{3.5}
For any $k\in \Char$, $(R_k,\chi_k)$ is a graded root.
\end{props}

\begin{proof}
Clearly, $\#\calv(R_k)$ is infinite. Notice
that the validity of Definition~\ref{2.1}(1) part (b)  guarantees that
$R_k$ has no
closed cycles. Indeed, in the presence of a cycle $v_1,v_2,
\ldots, $ $v_n,v_1$ in $R_k$, $\chi_k$ restricted to the set
$\{v_1,\ldots, v_n\}$ would have a local minimum which would
contradict Definition~\ref{2.1}(1) part (b). Hence, it is enough to verify
properties Definition~\ref{2.1}(1) parts (a)--(c). Property (a) is
clear. With the notation of (b), assume that (b) is not true for some
$u,v,w$, ie
$\chi_k(u)< \min\{\chi_k(v),\chi_k(w)\}$. This would imply
$C_u\subset C_v\cap C_w$. But this is impossible: two different
connected components of a space cannot simultaneously contain a
non-empty connected subset. The first two statements of (c) follow
from the fact that $B$ is negative definite. Finally, we have to
show that  $\overline{L}_{k,\leq n}$ is connected for $n$ sufficiently
large. Set
$$\caly:=\{x\in L:\, \chi_k(x)\leq \chi_k(x+\epsilon_jb_j)
  \ \mbox{for all}\ j\in \calj,\ \mbox{and}\ \epsilon_j=\pm 1\}.$$
Then, using the definition of $\chi_k$,
$\caly$ is
$$\{x\in L:\, -\chi_k(-b_j)\leq (x,b_j)\leq \chi_k(b_j)
  \ \mbox{for all}\ j\in\calj\}.$$
Hence $\caly$ is finite, and
obviously non-empty since $\chi_k$ has a global minimum (which is
in $\caly$). Fix $y_*\in\caly$, and for each
$y\in\caly\setminus\{y_*\}$ fix a path which connects $y_*$ and
$y$. By this we mean a sequence of elements $x_1,x_2,\ldots,x_t\in
L$ so that $x_1=y_*$, $x_t=y$, and
$x_{l+1}=x_l+\epsilon_{l}b_{j(l)}$ for some
$\epsilon_{l}\in\{+1,-1\}$. Let $n_0$ be the maximum of all the
values of type $\chi_k(x_l)$, where $x_l$ is an element on one of
the above paths connecting $y_*$ with some $y$. Then
$\overline{L}_{k,\leq n}$ is connected provided that $n\geq n_0$.
\end{proof}

Hence, for any $k\in \Char$ we get a
graded root  $(R_k,\chi_k)$. Some of these graded  roots  are not
very different. Indeed,  assume that $k$ and $k'$ determine the
same $\spin^c$ structure, hence $k'=k+2l$ for some $l\in L$. Then
$\chi_{k'}(x-l)=\chi_k(x)-\chi_k(l)$ for any $x\in L$. This means
that the transformation $x\mapsto x':=x-l$ realizes an
identification of $ \overline{L}_{k,\leq n}$ and $\overline{L}_{k', \leq
n-\chi_k(l)}$. Hence, we get:

\begin{props}\label{3.7}
If $k'=k+2l$ for some $l\in L$, then
$$(R_{k'},\chi_{k'})=(R_k,\chi_k)[-\chi_k(l)].$$
\end{props}

In fact,  there is an easy way to choose a  graded root from the
multitude $\{(R_k,\chi_k)\}_{k\in[k]}$.  For any $k\in \Char$ we define
\begin{equation}
m_k:=\frac{ k^2-\max_{k'\in[k]}(k')^2 }{8}\leq 0.\label{eq1}
\end{equation}
Clearly, $m_k$ is an integer (see the proof of the next lemma).  Set  also
$M_{[k]}:=\{k\in [k]:  m_k=0\}$.

\begin{lems}\label{3.8}
Fix a   $\spin^c$ structure  $[k]$.
Then $k_0\in M_{[k]}$ if and only if $-\chi_{k_0}(l)\leq 0$ for
any $l\in L$. Moreover, if $k_0$ and $k_0+2l\in M_{[k]}$, then
$-\chi_{k_0}(l)=0$.
\end{lems}

\begin{proof}
Write $(k_0+2l,k_0+2l)=(k_0,k_0)-8\chi_{k_0}(l)$.
\end{proof}

\begin{hh}\label{3.9} In particular, using Proposition~\ref{3.7} and
Lemma~\ref{3.8}, for
any fixed $\spin^c$ structure $[k]$, any choice of $k_0\in M_{[k]}$
provides  the same graded root $(R_{k_0},\chi_{k_0})$, which will
be denoted by $(R_{[k]},\chi_{[k]})$. Moreover, for any $k\in[k]$
\begin{equation}
(R_k,\chi_k)=(R_{[k]},\chi_{[k]})[m_k],\label{eq2}
\end{equation}
where $m_k$ is defined in \eqref{eq1}.
The notation is compatible with Proposition~\ref{2.7}:
\begin{equation}
m_k=\min\chi_k\label{eq3}
\end{equation}
Indeed, if $k_0\in M_{[k]}$, then $m_{k_0}=0$ from  \eqref{eq1}, and
$\min \chi_{k_0}=0$ from Proposition~\ref{3.8}.  Then \eqref{eq3} follows from
\eqref{eq2}. The relations  \eqref{eq1} and \eqref{eq3} can be summarized in
\begin{equation}
k^2-8\min \chi_k=\max_{k'\in[k]}(k')^2.\label{eq4}
\end{equation}\end{hh}

Clearly, many different plumbing graphs can provide the
same 3--manifold $M$. But all these plumbing graphs can be
connected by each other by a finite sequence of blowups/downs
$(-1)$--vertices, whose number of incident edges is $\leq 2$. (This
fact follows from the existence of an unique minimal resolution
graph of a normal surface singularity.)

The next result shows that the set
$\{(R_{[k]},\chi_{[k]})\}_{[k]}$ depends only on $M$.

\begin{props}
\label{3.10}
The set  $(R_{[k]},\chi_{[k]})$
(where $[k]$ runs over the $\spin^c$ structures of $M$) depends
only on $M$ and is independent of the choice of the (negative
definite) plumbing graph $\Gamma$ which provides $M$. (See also
\cite[Section~2.5]{OSzP}.)
\end{props}

\begin{proof}
By the above remark, it is enough to prove that the above set of
graded roots is not modified during a blowup. First we assume that
$\Gamma'$ is obtained from $\Gamma$ by ``blowing up a smooth point
of the Riemann surface'' which   in the plumbing construction
corresponds to a fixed index $j_0\in\calj$.  More precisely,
$\Gamma'$ denotes a graph with one more vertex and one more edge
than $\Gamma$: we glue to the vertex $j_0$ by the new edge the new
vertex with decoration $-1$, the decoration of $b_{j_0}$ is
modified from $e_{j_0}$ into $e_{j_0}-1$,  and we keep all the
other decorations. We will use the notation $L(\Gamma),\
L(\Gamma'),\ L'(\Gamma), \ L'(\Gamma'), \ B, \ B'$ for the
corresponding invariants.  Set $b_\text{new}$ for the new base element
in $L(\Gamma')$ (with $B'(b_\text{new},b_\text{new})=-1$).  The following
facts can be verified:

\begin{itemize}
\item Consider the maps $\pi_*\co L(\Gamma')\to
L(\Gamma)$ defined by
$$\pi_*\Bigl(\sum x_jb_j+x_\text{new}b_\text{new}\Bigr)=\sum x_jb_j,$$
and $\pi^*\co L(\Gamma)\to L(\Gamma') $  defined by
$$\pi^*\Bigl(\sum x_jb_j\Bigr)=\sum x_jb_j +x_{j_0}b_\text{new}.$$
Then $B'(\pi^*x,x')=B(x,\pi_*x')$.  This shows that
$B'(\pi^*x,\pi^*y)=B(x,y)$ and $B'(\pi^*x,b_\text{new})=0$ for any
$x,y\in L(\Gamma)$.

\item Identify (for both graphs) $L'$ with a
sub-lattice of $L_\Q$ (see Section~\ref{3.1}). Then
consider the (nonlinear) map $c\co L'(\Gamma)\to L'(\Gamma')$
defined by
$$\sum r_jb_j\mapsto \sum r_jb_j+(r_{j_0}+1)b_\text{new}$$
where $r_j\in \Q$, or
equivalently, $c(l'):=\pi^*(l')+b_\text{new}$.
Then $c(\Char(\Gamma))\subset \Char(\Gamma')$ and furthermore $c$  induces an
isomorphism between the orbit spaces $\Char(\Gamma)/2L(\Gamma)$ and $
\Char(\Gamma')/2L(\Gamma')$.

\item Consider $k\in \Char(\Gamma)$ and write $k':=c(k)\in
\Char(\Gamma')$.  Then for any $x\in L(\Gamma)$ one has
$\chi_{k}(x)=\chi_{k'} (\pi^*x)$.  In particular,  one has a
well-defined injection
$$\pi^*\co \{x\in L(\Gamma):\, \chi_{k}(x)\leq n\}
  \hookrightarrow \{y\in L(\Gamma'):\, \chi_{k'}(y)\leq n\}.$$

\item For any $z\in L(\Gamma')$ write $z$ in the
form $\pi^*\pi_*z+ab_\text{new}$ for some $a\in \Z$. Then
$\chi_{k'}(z)=\chi_{k'}(\pi^*\pi_*z)+\chi_{k'}(ab_\text{new})$. On the
other hand, for any $a'\in\Z$ (with $|a'|\leq |a|$), one has
$\chi_{k'}(a'b_\text{new})=a'(a'+1)/2\geq 0$. In particular, if
$\chi_{k'}(z)=m$, then $z$ and $\pi^*\pi_*z$ are in the same
connected component of $ \{y\in L(\Gamma'):\, \chi_{k'}(y)\leq
m\}$.
\end{itemize}

All these facts together imply $R_k=R_{k'}$.
There is a similar verification in the case when one blows up
``an intersection point''
corresponding to two indices $i_0$ and $j_0$ with $B(b_{i_0},b_{j_0})=1$.
The details are left to the reader.
\end{proof}

Finally, we verify the compatibility of the two combinatorial
objects $\bH^+(\Gamma,[k])$ (see Section~\ref{bhplus}) and
$\bH(R_{[k]},\chi_{[k]})$ (see Section~\ref{2.6}). Recall that $s$ denotes
the number of vertices of $\Gamma$.

\begin{props}\label{hbhb}
For any $[k]\in \Spin^c(M)$ one has
$$\bH^+(\Gamma,[k])=\bH(R_{[k]},\chi_{[k]})
  \big[-\max_{k'\in[k]}\frac{(k')^2+s}{4}\big].$$
\end{props}

\proof  Fix an arbitrary $k_0\in M_{[k]}$. Then $[k]=k_0+2L$,
hence $\bH^+(\Gamma,[k])\subset \Hom(\Char,\calt^+_0)$ can be identified with a set of
maps $\{\phi_L\co L\to \calt^+_0\}$ which satisfy the following:
for any $l\in L$ and $b_j$ with $2n=(k_0+2l+b_j,b_j)$ one has
$U^n\cdot \phi_L(l+b_j)=\phi_L(l)$ if $n\geq 0$, or $\phi_L(l+b_j)=U^{-n}\cdot \phi_L(l)$
if $n\leq 0$. Notice that $\chi_{k_0}(l+b_j)-\chi_{k_0}(l)=-n$; hence the above property
is equivalent with $(*)$ in Section~\ref{2.6}, showing the compatibility
of the two sets of restrictions. Still, one has to verify two more facts.

In order to formalize the first,  let us consider the natural map
$\theta\co L\to \calv(R_{k_0})$, where we associate with any $l\in L$ the
component of $\overline{L}_{k_0,\leq \chi_{k_0}(l)}$ containing $l$. This
induces a map $\Hom(\calv(R_{k_0}),\calt^+_0)\to \Hom(L,\calt^+_0)$
by composition $\phi\mapsto \phi_L:=\phi\circ \theta$ and, by the
above compatibility verifications, a map $\theta^*\co \bH(R_{k_0})\to
\bH^+(\Gamma,[k])$.  Since (in general) $\theta$ is not onto, one
has to verify that any $\phi_L$ has an extension to a $\phi$ with
$\theta^*(\phi)=\phi_L$, and this extension is unique.

Assume that $v\in \calv(R_{k_0})\setminus \im\theta$ corresponds
to a component $C_v$ in $\overline{L}_{k_0,\leq n}$. Set
$m_v:=\max(\chi_{k_0}(L\cap C_v)$, and take some $l_v\in L\cap C_v$ with
$\chi_{k_0}(l_v)=m_v$. Then for any $\phi\in \bH(R_{k_0})$, the value
$\phi(v)$ is uniquely determined by $\phi(v)= U^{n-m_v}\cdot \phi_L(l_v)$.

Finally, we have to fit the gradings. With the obvious notation,
for any $k=k_0+2l$, one has
$$\deg(\phi(k))=\frac{(k_0+2l)^2+s}{4}+d_{\bH^+}=-2\chi_{k_0}(l)+d_{\bH}$$
or, equivalently,
$$d_{\bH}-d_{\bH^+}=\frac{k_0^2+s}{4}=\max_{k'\in[k]}\frac{(k')^2+s}{4}.\eqno{\qed}$$

\section{Distinguished characteristic elements}

\noindent There is a more subtle way to choose a special characteristic
 element from each fixed  orbit $[k]$, a fact which will be
crucial in the next discussions. The goal of the present section
is its definition.

\begin{h}\textbf{Definitions: The canonical graded root}\qua
\label{4.1}
Recall that the {\em canonical} characteristic element  $K\in \Char$  is
defined by the  equations $K(b_j)=-e_j-2$
(or equivalently, by $\chi_K(b_j)=1$)    for any
$j\in \calj$. We denote by $(R_\text{can},\chi_\text{can})$ the graded root
associated with $K$.
In order to simplify the notation,  we also write $\chi_K=\chi$.

Sometimes it is preferable (motivated by the symplectic geometry)
to define the ``canonical'' $\spin^c$ structure $\sigma_\text{can}$ via
$c_1(\sigma_\text{can})=
[-K]$, ie identifying it by the anti-canonical bundle. Our choice above is motivated by the
Riemann--Roch formulas, and the connection with singularity theory.
 But, in fact, from the point of view of the results of the present paper,
the two choices are completely equivalent. Indeed, there is a
natural involution $x\mapsto -x$ on $L$,  $L'$ and $H$. This is
compatible with the natural involution on $\Spin^c(M)$. All the
formulas in this paper are  stable with respect to these
involutions (eg, $\min \chi_K=\min\chi_{-K}$): the orbits $[K]$
and $[-K]$ share the same properties.
\end{h}

\begin{exss}\label{4.2}
If $M$ is the link of a normal surface singularity, and
$\widetilde{X}$ is a resolution of $(X,0)$, then $K$ is the first Chern class
of the canonical bundle of the complex structure of $\widetilde{X}$.
$b_j\in L$ denotes  the fundamental class of the
irreducible exceptional divisors $E_j$, hence
$x=\sum n_jb_j\in L$ can be identified with the cycle $Z=\sum n_jE_j$
supported by the exceptional divisor.
 The linear equations defining $K$ are
the adjunction formulas for the canonical line bundle.  Hence, by
Riemann--Roch, $\chi(x)$ is exactly the Euler-characteristic of the
sheaf-cohomology of ${\mathcal O}_Z$.\end{exss}

\begin{h}\textbf{Notation}\qua
\label{n1}
We embed $L'$ into $L_\Q$ as in  Section~\ref{3.1}.
We set
$$S_\Q:=\{x\in L_\Q:\, (x,b_j)\leq 0 \ \mbox{for every $j\in
\calj$}\}$$ and $S=S_\Q\cap L$.
Since $B$ is negative definite, if
$x\in S_\Q$, then $x\geq 0$, see Theorem~\ref{4.3} part (a).
\end{h}

\begin{h}\label{luj2}
Fix a characteristic element $k$, ie fix an
$l'\in L'$ with $k=K+2l'$.
Then the  function $\chi_k\co L\to \Z$
is equivalent to the rational valued extension of $\chi=\chi_K$
restricted to the sub-lattice $l'+L$ of $L_\Q$. Indeed, extend $\chi$ by the
same formula $\chi\co L_\Q\to \Q$, $\chi(x):=-(K+x,x/2)$. Then, for any
$x\in L$ one has
 \begin{equation*}
\chi_k(x)=\chi(l'+x)-\chi(l').
\end{equation*}
In fact, the equivalence class $[k]=K+2l'+2L$ can be identified
with this sub-lattice   $l'+L$ in $L_\Q$. We wish to choose a
distinguished element of $l'+L$ and list some of its
properties.\end{h}

\begin{lems}\label{luj}
For any fixed $[k]=K+2(l'+L)$,
the intersection $(l'+L)\cap S_\Q$ in $L_\Q$ admits a unique
minimal element $l'_{[k]}$. (Here minimality is considered with
respect to the ordering $\leq $ in $L_\Q$.)
\end{lems}
\begin{proof}
This is similar to Artin's proof of the
existence of the fundamental cycle \cite{Artin62,Artin66}.
Assume that $l'+l_i\in S_\Q$, $l_i\in L$ for $i=1,2$. Set
$l:=\min\{l_1,l_2\}$, and write $x_i=l_i-l\geq 0$.  We show that
$l'+l\in S_\Q$ as well.  Consider any $j\in\calj$. Since
$|x_1|\cap|x_2|=\emptyset$, we can assume that either $b_j\not\in
|x_1|$ or  $b_j\not\in |x_2|$. If $b_j\not\in |x_1|$ then
$(l'+l,b_j)=(l'+l_1,b_j)-(x_1,b_j)\leq 0$.
\end{proof}

\begin{defs}\label{defuj}
In any fixed class $[k]$ we fix the distinguished
representative $k_r=K+2l'_{[k]}\in[k]$.
\end{defs}


\begin{h}\textbf{Remarks}\qua
\label{ineq}
(a)\qua Since $l'_{[K]}=0$, the distinguished representative
in $[K]$ is $K$.

(b)\qua If $e_j\leq -2$ for any $j$ (eg, if $\Gamma$ is a {\em
minimal} resolution graph of a normal surface singularity) then
$K(x)\geq 0$ for any $x\geq 0$. This fails to be true  for
arbitrary plumbing graphs.  Still, we have $K(x)-(x,x) \geq  0 $
for any $x \geq 0$. This is a special case of the following
general result.
\end{h}

\begin{props}\label{x1}
For any fixed $[k]\in \Spin^c(M)$ and
the representative $k_r\in[k]$ one  has:

{\rm(a)}\qua $(l'_{[k]},b_j)\geq b_j^2+1$ for any $j\in\calj$;

{\rm(b)}\qua $ k_r(x)\geq x^2$ for any $x\geq 0$, $ x\in L$;

{\rm(c)}\qua $\chi_{k_r}(-x)\geq 0$ for any $x\geq 0$, $x\in L$.
\end{props}

\begin{proof} (a) Assume that $(l'_{[k]},b_{j_0})\leq b_{j_0}^2$ for  some $j_0\in\calj$.
Then by a computation $l'_{[k]}-b_{j_0}\in (l'_{[k]}+L)\cap S_\Q$, a fact
which contradicts the minimality of $l'_{[k]}$, see Lemma~\ref{luj}.

Part (b) follows similarly. First notice that by (a)
 $k_r(b_j)= 2(l'_{[k]},b_j)-b_j^2-2\geq b_j^2$
for any $j\in\calj$.
Therefore, if there exists an effective cycle $x>0$ with
$k_r(x)<x^2$,  then there exists a minimal one with this property.
This means that that minimal $x>0$ satisfies
$$x^2-k_r(b_j)>k_r(x)-k_r(b_j)\geq (x-b_j)^2 \ \ \ \mbox{for every $b_j\in |x|$.}$$
Hence $0\geq (l'_{[k]}-x,b_j)$ for every $b_j\in|x|$. On the other
hand, if $b_j\not\in|x|$, then $0\geq (l'_{[k]},b_j) \geq
(l'_{[k]}-x,b_j)$. These two inequalities show that for this $x>0$
one has $l'_{[k]}-x\in (l'_{[k]} +L)\cap S_\Q$, which contradicts
the minimality of $l'_{[k]}$. For (c) use (b) and the identity
$2\chi_{k_r}(-x)=k_r(x)-x^2$.
\end{proof}

\section{First  properties of the graded roots
$(R_k,\chi_k)$}\label{sec6}

\noindent  In this section we focus on the distinguished
representatives $k_r=K+2l'_{[k]}$ of the orbits $[k]$.
(Then all their properties can be transformed into
similar properties of arbitrary characteristic elements via
Proposition~\ref{3.7}.)

 For any class $[k]$, we define
$S_{[k]}:=\{x\in L:\, (x+l'_{[k]},b_j)\leq 0 \  \mbox{for every
$j\in\calj$}\}.$ Clearly, $S_{[K]}=S$.

\begin{thms}\label{4.3}
{\rm(a)}\qua For any $[k]$, if $x\in S_{[k]}$ then $x\geq 0$.

{\rm(b)}\qua For any class $[k]$,  consider the distinguished representative
$k=k_r=K+2l'_{[k]}$.
Let $\overline{L}^s_{k,\leq n}$ be the subcomplex of $
\overline{L}_{k,\leq n}$ whose 0--skeleton consists of  cycles $x\in
L^s_{k,\leq n}:=L_{k,\leq n}\cap S_{[k]}$, and the 1--cells are exactly
those 1--cells of
$\overline{L}_{k,\leq n}$ which have their endpoints
in $L^s_{k,\leq n}$. Then the natural inclusion
$\overline{L}^s_{k,\leq n}\to \overline{L}_{k,\leq n}$
induces a surjection
$\pi_0(\overline{L}^s_{k,\leq n})\to \pi_0(\overline{L}_{k,\leq n})$ for any $n$.

{\rm(c)}\qua Assume $k=K$. Then the component $C_0$ in $\overline{L}_{K,\leq 0}$,
 which contains the zero cycle,  contains no non-zero effective cycle.
In fact, if ${\mathcal X}$ denotes  the set of non-zero cycles
in the component $C_0$, then for any $x\in{\mathcal X}$ one has $x<0$
and  $\chi(x)=0$.

{\rm(d)}\qua  For $k=K$, $\overline{L}_{K,\leq n}$ is connected whenever $n\geq
1$.
\end{thms}

\begin{proof}  (a)\qua Notice that $(l'_{[k]}+L)\cap S_\Q=l'_{[k]}+S_{[k]}$,
hence (a)  follows from the minimality of $l'_{[k]}$, see Lemma~\ref{luj}.
A direct argument goes as follows.
Assume the contrary: $x=x_1-x_2$ with $x_1\geq 0$, $x_2>0$ and $|x_1|\cap|x_2|=
\emptyset$.
If $b_j\in|x_2|$, then $(b_j,x_1)\geq 0$, hence by the definition of $S_{[k]}$ one gets
 $(l'_{[k]}-x_2,b_j)\leq -(x_1,b_j)\leq 0$.
If $b_j\not\in |x_2|$ then $(x_2,b_j)\geq 0\geq  (l'_{[k]},b_j)$.
Hence $l'_{[k]}-x_2\in S_\Q$, which contradicts the minimality of
$l'_{[k]}$, see Lemma~\ref{luj}.

(b)\qua Assume that $x\in L\setminus S_{[k]}$.
  Then for some
$j\in\calj$, $(x+l'_{[k]},b_j)\geq 1$. This is equivalent with
$\chi_k(x+b_j)\leq \chi_k(x)$.
If $x+b_j\in S_{[k]}$ then we stop, otherwise we repeat the same algorithm for
$x+b_j$ instead of $x$. In this way we construct an increasing sequence along which
$\chi_k$ is decreasing. Since $B$ is negative definite, this procedure
must stop.

(c)\qua We prove that any $x\in \calx$ satisfies
\begin{enumerate}
\item[(i)] $x<0$ and
\item[(ii)] $\chi(x)=0$.
\end{enumerate}
First notice that $b_j\not\in
L_{K,\leq 0}$ and $-b_j\in L_{K,\leq 0}$ if and only
if $\chi(-b_j)=0$, or equivalently, $e_j=-1$.

Since $\calx$ represents the (non-zero) cycles of the connected component
$C_0$, it is enough to verify the following inductive step.
If $x\in\calx$ satisfies (i) and (ii), and $x'=x\pm b_j\in\calx$, then
$x'$ satisfies (i) and (ii) as well. In order to verify (i) for $x'$,
we have to verify that the situation $x'=x+b_j$ with $b_j\not\in|x|$
is not possible. Indeed, in this case $\chi(x')=\chi(x)+\chi(b_j)-(x,b_j)=
1+(-x,b_j)\geq 1$ which would contradict $x'\in L_{K,\leq 0}$. Now, using
(i) and the fact $\chi(x')\leq 0$, Proposition~\ref{x1} part (c) implies
$\chi(x')=0$.

(d)\qua Via (a) and (b), it is enough to  show that any $x_0\geq 0$,
$x_0\in L_{K,\leq n}$,  can be connected to the zero cycle by a path
inside of $\overline{L}_{K,\leq n}$. For this, it is enough to verify
the following two facts:

(i)\qua For any $x>0$ with $\chi(x)>0$, there is a $b_j\in |x|$
such that $\chi(x-b_j)\leq \chi(x)$.

Indeed, if this were not the case, then we would have $(b_j,x-b_j)>1$
for any $b_j\in|x|$. This is equivalent with $(b_j,x)\geq -K(b_j)$,
which implies $K(x)+x^2\geq 0$, or $\chi(x)\leq 0$, a contradiction.

(ii)\qua For any $x>0$ with $\chi(x)\leq 0$, there exists $b_j\in |x|$
such that  $\chi(x-b_j)\leq 1$.

To prove this, denote $\chi(x)$ by $c\leq 0$,
write $x=\sum_jn_jb_j$ and  set $n:=\sum_jn_j$. Clearly $n\geq 2$ since
$\chi(b_j)=1$ for any $j$. Now,
assume that $\chi(x-b_j)>1$ for all $b_j\in |x|$. Then,
in the identity $-(b_j,x-b_j)=\chi(x)-\chi(x-b_j)-1$ the right
hand side is $<c-2$. This means that $-(b_j,x)<c-2-e_j=c+K(b_j)$
for any $b_j\in|x|$. This implies that $K(x)+x^2+cn>0$, or $c(n-2)>0$,
which is a contradiction.

Then the construction of the wanted path connecting $x_0$ with
$0$ is clear: in the path any $x$ is succeeded by some $x-b_j$ according
to (i) or (ii).
\end{proof}

\begin{h}\textbf{Characterization of rational singularities}\qua
\label{4.4}
Recall
that a normal surface singularity is called {\em rational} if its geometric
genus $p_g=0$. It is easy to verify that any resolution of a rational
singularity is automatically good (ie any resolution graph
can be considered as a plumbing graph), and the link $M$ is a rational homology sphere.
Artin characterized the rational singularities topologically in terms of
{\em  any} fixed resolution graph $\Gamma$ \cite{Artin62,Artin66}. Namely,
\begin{equation*}
\mbox{$p_g=0$ if and only if $\chi(x)\geq 1$ for any $x>0$, $x\in L$.}
\tag{1}\end{equation*}
A connected, negative definite plumbing graph with this property is
called a {\em rational graph}.

For example, any plumbing graph with $-e_j\geq \delta_j$ for every
$j\in\calj$ is rational.  This can be verified as follows: if
$-e_j\geq \delta_j$ then the Artin's fundamental cycle $x_\text{min}$
(see Remark~\ref{artin}) is exactly $\sum_jb_j$, and by an easy
computation one has  $\chi(x_\text{min})=1$, hence the rationality
follows from Artin's criterion, see Remark~\ref{artin} and
\cite{Artin62,Artin66}.  On the other hand, not any rational graph
satisfies $-e_j\geq \delta_j$: consider eg the $-E_8$ graph.

The class of rational graphs is closed while taking subgraphs and
decreasing Euler numbers $e_j$.
\end{h}

\begin{thms}\label{4.5} Let $\Gamma$ be a connected,
negative definite plumbing graph whose plumbed three-manifold is a
rational homology sphere. Then the following facts are equivalent:

{\rm(a)}\qua $\Gamma$ is rational;

{\rm(a$'$)}\qua $\# \chi_\text{can}^{-1}(0)=1$ (where
$(R_\text{can},\chi_\text{can})$ is the canonical
graded root associated with $\Gamma$);

{\rm(b)}\qua $R_\text{can}=R_0$ (for the definition of $R_0$, or of any $R_m$,
see Example~\ref{2.3}(1));

{\rm(c)}\qua $R_\text{can}=R_m$ for some $m\in\Z$;

{\rm(d)}\qua For all characteristic elements  $k\in \Char$, $R_k=R_{m_k}$ for some $m_k\in \Z$;

{\rm(Hb)}\qua $\bH(R_\text{can})=\calt^+_0$;

{\rm(Hc)}\qua $\bH(R_\text{can})=\calt^+_m$ for some $m\in\Z$; or equivalently,
 $\bH_\text{red}(R_\text{can})=0$;

{\rm(Hd)}\qua For all $k\in \Char$, $\bH(R_k)=\calt^+_{m_k}$ for some $m_k\in \Z$;
or equivalently,  $\bigoplus_{[k]}\bH_\text{red}(R_{[k]})=0$.

Moreover, if $\Gamma$ is rational and $k=K+2l'$, then
$$m_k=\min\chi_k=\min_{x\in L}\chi(l'+x)-\chi(l')=\chi(l'_{[k]})-\chi(l')=\chi_k(l'_{[k]}-l')\leq 0.$$
In particular, if $\Gamma$ is rational and $k_r=K+2l'_{[k]}$, then
$\min \chi_{k_r}=0$ (see also Proposition~\ref{3.7} and Section~\ref{3.9}).
\end{thms}

It is instructive to compare (a') with the property
$\#\chi_\text{can}^{-1}(n)=1$,  valid  for any graph $\Gamma$ and integer
$n\geq 1$ (see Theorem~\ref{4.3} part (d)).

\begin{proof} (a) $\Rightarrow$ (b) follows from (1) in Section~\ref{4.4} and
Theorem~\ref{4.3} (since $\overline{L}^s_{K,\leq 0}$ consists of only
one vertex, namely $x=0$, and
$\overline{L}^s_{K,\leq n}=\emptyset$ for $n<0$). (b) $\Rightarrow$
(c) $\Rightarrow$ (Hc) and (b) $\Rightarrow$ (Hb) $\Rightarrow$ (Hc) are clear.
(Hc) $\Rightarrow$ (c) follows from Proposition~\ref{2.7}. Next, we verify
(c) $\Rightarrow$ (a') $\Rightarrow$ (a). Indeed, if $R_\text{can}=R_m$  for some
$m$, then clearly $m\leq \chi(0)=0$
 and $\chi_\text{can}^{-1}(0)$ has only one element. This clearly corresponds
to the connected component $C_0$ of $\overline{L}_{K,\leq 0}$ which contains the
zero cycle. Hence $L_{K,\leq 0}=C_0$,
and  by Theorem~\ref{4.3} part (c)  $\chi(x)>0$ for any $x>0$.

Finally, it is clear that (d) $\Rightarrow$ (Hd) $\Rightarrow$ (Hc). Hence
it remains to show that (a) $\Rightarrow$ (d).  Using
Proposition~\ref{3.7}, it is enough to consider the
characteristic elements of the form $k=K+2l'_{[k]}$. Fix such a $k$.

Let $L^e_{k,\leq n}=L_{k,\leq n}\cap \{\mbox{effective cycles}\}$,
and define $\overline{L}^e_{k,\leq n}$  similarly as $\overline{L}^s_{k,\leq
n}$ in Theorem~\ref{4.3} part (b).  Then Theorem~\ref{4.3} parts (a)
and (b) show that
$\pi_0(\overline{L}^e_{k,\leq n})\to \pi_0(\overline{L}_{k,\leq n})$ is onto
for any $n$. Hence it is enough to show that any $x\in L^e_{k,\leq
n}$ can be connected with the zero cycle in $\overline{L}^e_{k,\leq n}$
(a fact which automatically shows that $n\geq 0$ whenever
$L^e_{k,\leq n}$ is non-empty). More precisely, we will show that
for any $x>0$ there exists at least one $b_j\in|x|$ so that
$\chi_k(x-b_j)\leq \chi_k(x)$. Indeed, assume that this is not
true for some $b_j$. Then this means
$(x+l'_{[k]},b_j)>1+b_j^2$ for every $b_j\in|x|$. Since
$(l'_{[k]},b_j)\leq 0$, we get $(K+x,b_j)\geq 0$ for every
$b_j\in|x|$. This implies $K(x)+x^2\geq 0$, or  $\chi(x)\leq 0$,
 which contradicts the fact that $\Gamma$ is rational.

This, in particular, also shows that $\chi_{k_r}(x)\geq 0$ for any
$x\geq 0$. From the proof of Theorem~\ref{4.3} part (b)
 we get that $\chi_{k_r}(x)\geq 0$ for
any $x$. In other words, $\chi(l'_{[k]}+x)\geq \chi(l'_{[k]})$ for
any $x\in L$. This and Section~\ref{luj2} also imply the last statement
about $m_k$.
\end{proof}

\begin{h}\textbf{Characterization of weakly elliptic singularities via
$(R_\text{can},\chi_\text{can})$ or $\bH(R_\text{can})$}\qua
\label{4.6}
A normal
surface singularity is called \emph{weakly elliptic} if any of its
resolution graphs  are weakly elliptic. A resolution graph
$\Gamma$ is weakly elliptic if $\min_{x>0} \chi(x)=0$ (see
\cite{Laufer77}). (This definition is valid even if the link is
not a rational homology sphere, but in order to be consistent with
all the above notation, we will assume this fact as well.)\end{h}

\begin{props}\label{4.7}
Let $\Gamma$ be a connected,
negative definite plumbing graph whose plumbed three-manifold is a
rational homology sphere. Then the following facts are equivalent:

{\rm(a)}\qua $\Gamma$ is weakly elliptic.

{\rm(b)}\qua $R_\text{can}=R(\{n_i\},\{n_{ij}\})$ for some index set $I$, $\# I=l+1
\geq 2$, and $n_i=0$ for any $i\in I$, and $n_{ij}=1$ for any $i\not=j$.

{\rm(Hb)}\qua $\bH(R_\text{can})=\calt^+_0\oplus \big(\,
\calt_0(1)\,\big)^{\oplus l}$ for some $l\geq 1$.
\end{props}

\begin{proof} (a) $\Leftrightarrow$ (b) follows from
Theorem~\ref{4.3}, and  (b) $\Leftrightarrow$ (Hb) from
Theorem~\ref{4.3}(d) and Proposition~\ref{2.7}. \end{proof}

\begin{h}\textbf{Remarks}\qua
\label{re}
(a)\qua The results Theorem~\ref{4.3} part (d),
Theorem~\ref{4.5} and Proposition~\ref{4.7} can also be
interpreted as follows: The grading $\chi_\text{can}$ of
$R_\text{can}$ satisfies $\min \chi_\text{can}\geq 0$ (or, equivalently
$\min \chi_\text{can}=0$)
if and only if $\Gamma$
is either rational or weakly elliptic. In this situation, if
$\bH_\text{red}(R_\text{can})=0$ then $\Gamma$ is rational, otherwise it is weakly
elliptic.

(b)\qua Let $k_r=K+2l'_{[k]}$ as above. Then $\min\chi_{k_r}\geq \min\chi$.
Indeed, by Theorem~\ref{4.3}(a)--(b), $\min\chi_{k_r} = \chi_{k_r}(x)$ for some
effective cycle $x\geq 0$. But for such a cycle $ \chi_{k_r}(x)=\chi(x)-(l'_{[k]},x)\geq \chi(x)$.

(c)\qua In particular, if $\Gamma$ is rational or weakly elliptic,
then for any $\spin^c$--structure $[k]$ one has $\min\chi_{k_r}=0$.
Indeed, $0=\chi_{k_r}(0)\geq \min \chi_{k_r}\geq
\min\chi=0$.\end{h}

\begin{h}\textbf{Remark}\qua
\label{re2}
One has the following connections with singularity theory.
(We omit the details, since these facts are not exactly in the spirit
of the present paper.) One can verify that $l$ in Proposition~\ref{4.7}
is exactly the length of the elliptic  sequence in the sense of
S.S.-T. Yau, or, equivalently,  it is the length of Laufer's
sequence. (For their definitions and equivalence, see
\cite{Stevens}.) If the graph $\Gamma$  is minimally elliptic (in
the sense of Laufer \cite{Laufer77}) then $l=1$.  The opposite
statement is not necessarily true: one can find elliptic graphs
with $l=1$ which are not  numerically Gorenstein (ie $K$ is not
an integral cycle), hence which are not  minimally elliptic. For example,
the following graph $\Gamma$ has these properties:\end{h}

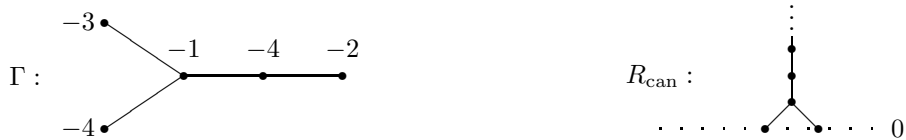
\begin{picture}(300,60)(0,10)\small
\linethickness{.5pt}
\put(70,40){\circle*{3}}
\put(100,40){\circle*{3}}
\put(130,40){\circle*{3}}
\put(40,60){\circle*{3}}
\put(40,20){\circle*{3}}
\put(40,60){\line(3,-2){30}}
\put(40,20){\line(3,2){30}}
\put(70,40){\line(1,0){60}}
\put(70,50){\makebox(0,0){$-1$}}
\put(100,50){\makebox(0,0){$-4$}}
\put(130,50){\makebox(0,0){$-2$}}
\put(30,60){\makebox(0,0){$-3$}}
\put(30,20){\makebox(0,0){$-4$}}
\put(10,40){\makebox(0,0){$\Gamma:$}}

\put(250,40){\makebox(0,0){$R_\text{can}:$}}

\put(300,40){\circle*{3}}
\put(300,50){\circle*{3}}
\put(300,30){\circle*{3}}
\put(290,20){\circle*{3}}
\put(310,20){\circle*{3}}
\put(300,30){\line(-1,-1){10}}
\put(300,30){\line(1,-1){10}}
\put(300,30){\line(0,1){25}}
\put(300,65){\makebox(0,0){$\vdots$}}
\put(340,20){\makebox(0,0){$0$}}
\qbezier[10](250,20)(290,20)(330,20)
\end{picture}

(Here the two minimal points of $R_\text{can}$  correspond to the zero cycle and
to Artin's fundamental cycle.)

\section{Generalities about computation sequences}\label{cs}

\begin{h}
The computation of the groups $\bH$ (see Section~\ref{sec9})
is based on the techniques of computation sequences (see, for example,
\cite{Laufer72,Laufer77,Ninv,Stevens}). These objects were
successfully used in the study of the (resolution of) normal
surface singularities (see \cite{Five} for more details).
Some of the next statements and proofs rhyme perfectly with some
of those computations, eg with the proof of the existence of
Artin's fundamental cycle, with Laufer's algorithm which provides
this fundamental cycle, or with the construction of Yau's elliptic
sequence.  Nevertheless, in order to be self-contained, and since
we also wish to treat the case of an arbitrary $\spin^c$--structure,
we will provide all the details.\end{h}

\begin{defs}
\label{gen1} Sequences
$x_0,x_1,\ldots, x_t\in L$ with
 $x_{l+1}=x_l+b_{j(l)}$ ($0\leq l<t$), where $j(l)$ is determined by some  principles
fixed in each individual case, are called {\em computation sequences} connecting $x_0$ and $x_t$.

In this section we will fix a $\spin^c$--structure $[k]$ and we will
fix its distinguished representative $k_r=K+2l'_{[k]}$. In
order to simplify the notation, in this section we write $k$ for
$k_r$.  Recall that $S_{[k]}=\{x\in L:\, (x+l'_{[k]},b_j)\leq 0\
\ \mbox{for every  $j\in\calj$}\}$.  
We will write $m_j(z)$ for the coefficients of a cycle $z=\sum
m_j(z)b_j$.\end{defs}

\begin{lems}
\label{genlem1}
For any $x\in L$, there exists a unique
minimal  $y=y(x)\in L$ with the properties  $x\leq y$ and $y\in
S_{[k]}$.
\end{lems}

\begin{proof}
First notice that (since $B$ is negative definite)  there exists a
cycle $y>0$ so that $(y,b_j)<0$ for every $j\in\calj$. Then  all
the coefficients of $y$ are strict positive (see also
Theorem~\ref{4.3}(a)). Hence, some integral multiple $ay$ of it will
satisfy $ay>x$ and $ay\in S_{[k]}$.  For the existence of a unique
minimal element with these properties it is enough to verify that
if $y_1,y_2\in S_{[k]}$ then $y:=\min\{y_1,y_2\}$ is also in
$S_{[k]}$. Indeed, for any $j$,  at least for one index
$i\in\{1,2\}$ one has $b_j\not\in |y_i-y|$. Then
$(y+l'_{[k]},b_j)=(y_i+l'_{[k]},b_j)-(y_i-y,b_j)\leq 0$.
\end{proof}

\begin{lems}{\rm(Generalized Laufer's algorithm \cite{Laufer72})}\qua
\label{genlem2}
Fix an $x\in L$.  Construct  a computation sequence $x_0,\ldots, x_t$ by
the following algorithm. Set $x_0=x$.
Assume that $x_l$ is already constructed.
 If $x_l\not\in S_{[k]}$, ie  $(x_l+l'_{[k]},b_j)>0$ for some index
$j$, then choose one of them
for $j(l)$, and write $x_{l+1}=x_l+b_{j(l)}$.
If $x_l\in S_{[k]}$, then stop and write  $t=l$.
Then $x_t$ is exactly $y(x)$  considered in Lemma~\ref{genlem1}.

Moreover,  this computation sequence satisfies
$\chi_k(x_{l+1})\leq \chi_k(x_l)$ for every $0\leq l<t$.
\end{lems}

\begin{proof}
We  will show by induction that $x_l\leq y(x)$ for any $0\leq l<t$.  Indeed, assume that
$x_l\leq y(x)$, and $x_{l+1}=x_l+b_{j}$.  Then we have to verify that
$m_j(x_l)<m_j(y(x))$.  Assume that this is not true, ie  $b_j\not\in |y(x)-x_l|$. Then
$(x_l+l'_{[k]},b_j)=(y(x)+l'_{[k]},b_j)-(y(x)-x_l,b_j)\leq 0$, a contradiction.

On the other hand, if each $x_l\leq y(x)$, then the algorithm should stop with
some $x_t\in S_{[k]}$, and $x_t\leq y(x)$. Then the minimality of $y(x)$ guarantees
that $x_t=y(x)$.

The inequality follows from
$\chi_k(x_{l+1})-\chi_k(x_l) =1-(x_l+l'_{[k]},b_j)\leq 0$.
\end{proof}

\begin{h}\textbf{Remark}\qua\label{artin}
(a)\qua If we start with $x=b_j$ for some $j$,  and $k=K$, then
$y(x)$ is exactly Artin's fundamental cycle $x_\text{min}$. It is the
minimal strict effective cycle with $(x_\text{min},b_j)\leq 0$ for
every $j$. Since $x_\text{min}\geq \sum_jb_j$, hence also $x_\text{min}\geq
b_j$ for every $j\in \calj$, the proof of Lemma~\ref{genlem2} also shows
that starting with $x_0=b_j$ we get $x_\text{min}$ independently of the
choice of the index $j\in\calj$. In fact, there is a computation
sequence connecting $x_0=b_j$ and  $x_\text{min}$ such that
$x_{s-1}=\sum_jb_j$  and $\chi(x_l)=1$ for any $0\leq j <s$, where
$s=\#\calj$.

It is known that $\Gamma$ is rational if and only if
$\chi(x_\text{min})=1$ (see \cite{Artin62,Artin66}), and $\Gamma$ is weakly
elliptic if and only if $\chi(x_\text{min})=0$ (see \cite{Laufer77}). In
general, $\chi(x_\text{min})\leq 1$.

(b)\qua Assume that above $\Gamma$ is  a rational graph,
and $k=K$ (hence $l'_{[k]}=0$ and $\chi_k=\chi$). Consider the
computation sequence from Lemma~\ref{genlem2} connecting $x_0:=b_j$ and
$x_\text{min}$.  It satisfies $\chi(x_{l+1})\leq \chi(x_l)$ for any
$l$. On the other hand $\chi(b_j)=1$ and $\chi(x_\text{min})= 1$, by
the rationality of $\Gamma$. Hence $\chi(x_{l+1})=\chi(x_l)$ for
any $l$; in other words, $(x_l,b_{j(l)})=1$ for any $l$
 (a fact first noticed by Laufer \cite{Laufer72}).\end{h}

In the next paragraphs we wish to generalize  the above  lemmas
\ref{genlem1} and \ref{genlem2}.
We prefer to write $\calj =\{0\}\cup \calj^* $,  and hence  to  distinguish
a base element $b_0$.

\begin{lems}\label{genlem3}
For any integer $i\geq 0$, there exists a
unique cycle $x(i)\in L$ with the following properties:

{\rm(a)}\qua $m_0(x(i))=i$;

{\rm(b)}\qua $(x(i)+l'_{[k]}, b_j)\leq 0$ for any $j\in\calj^*$;

{\rm(c)}\qua $x(i)$ is minimal with properties  {\rm(a)}--{\rm(b)}.

Moreover, the cycle $x(i)$ satisfies $x(i)\geq 0$.
\end{lems}

\begin{proof}
First we verify that there exists at least one $x\in L$ satisfying
(a) and (b). By Lemma~\ref{genlem1} there exists $\widetilde{x}\in S_{[k]}$ with
$\widetilde{x}\geq b_0$. Take an integer $a\geq 1$ such that
$(a-1)l'_{[k]}$ is an integral cycle and
$h:=m_0(a\widetilde{x}+(a-1)l'_{[k]})-i\geq 0$. Then define
$x:=a\widetilde{x}+ (a-1)l'_{[k]}-hb_0$. By a computation $m_0(x)=i$
and $(x+l'_{[k]},b_j)= a(\widetilde{x}+l'_{[k]},b_j)-h(b_0,b_j)\leq 0$
for $j\not=0$.

Next, one can verify that there exists a unique minimal element
with properties (a) and (b). The proof is similar to the proof of
Lemma~\ref{genlem1}, and it is left to the reader.

Finally we verify that $x(i)\geq 0$.   Write $x(i)$ in the form
$x_1-x_2$ with $x_1\geq 0$, $x_2\geq 0$, $|x_1|\cap
|x_2|=\emptyset$. Consider an index $j\in\calj^*$. If  $b_j\not\in
|x_1|$ then one has $(l'_{[k]}-x_2,b_j)\leq
(l'_{[k]}-x_2+x_1,b_j)= (x(i)+l'_{[k]},b_j)\leq 0$. Similarly, for
any $b_j\in |x_1|$ one has $(l'_{[k]}-x_2,b_j)\leq
(l'_{[k]},b_j)\leq 0$. Hence $l'_{[k]}-x_2\in (l'_{[k]}+L)\cap
S_\Q$, which implies $x_2=0$ via the minimality of $l'_{[k]}$, see
Lemma~\ref{luj}.
\end{proof}

\begin{lems}{\rm(The computation sequence connecting $x(i)$ and
$x(i+1)$)}\qua
\label{genlem4}
For any integer $i\geq 0$ consider a  computation sequence constructed as follows.
Set $x_0=x(i)$, $x_1=x(i)+b_0$. Assume that $x_l$ ($l\geq 1$)  is already constructed.
If $x_l$ does not satisfy property (b) in Lemma~\ref{genlem3}, then
there exists some $j\in\calj^*$
with $(x_l+l'_{[k]},b_j)>0$.  Then choose one of these indices
for $j(l)$, and write $x_{l+1}=x_l+b_{j(l)}$.
If $x_l$  satisfies property (b) in Lemma~\ref{genlem3}, then stop and
write  $t=l$.  Then $x_t$ is exactly $x(i+1)$.

Moreover,  this computation sequence satisfies
$\chi_k(x_{l+1})\leq \chi_k(x_l)$ for any $0<l<t$.

Corresponding to $l=0$ one has
$\chi_k(x(i)+b_0)-\chi_k(x(i))=1-(x(i)+l'_{[k]},b_0)$. In other
words, $\chi_k(x(i)+b_0)>\chi_k(x(i))$ if and only if $x(i)\in
S_{[k]}$.
\end{lems}

\begin{proof}
Repeat the arguments used in the proof of Lemma~\ref{genlem2}.
\end{proof}

Lemma~\ref{genlem4} has the following easy generalization (with the same proof):

\begin{lems}\label{genlem5}
Assume that $x\in L$ satisfies $m_0(x)=i$ and
$x\leq x(i)$ for some $i\geq 0$. Consider a similar   computation
sequence  as in Lemma~\ref{genlem4}. Namely, set $x_0=x$.  Assume that
$x_l$ is already constructed. If  for some $j\in\calj^* $ one has
$(x_l+l'_{[k]},b_j)>0$  then take $x_{l+1}=x_l+b_{j(l)}$, where
$j(l) $ is such an index $j$. If $x_l$  satisfies property (b) of
Lemma~\ref{genlem3}, then stop and write  $t=l$. Then $x_t$ is
exactly $x(i)$. Moreover,  this computation sequence satisfies
$\chi_k(x_{l+1})\leq \chi_k(x_l)$ for any $0\leq l<t$.
\end{lems}

Notice that, even if it is not explicitly emphasized in its notation,
the cycles $\{x(i)\}_{i\geq 0}$ depend on the choice of
the distinguished vertex and of the $\spin^c$--structure $[k]$.

\section{Almost-rational graphs,
$HF^+(-M,[k])=\bH^+(\Gamma,[k])$}\label{sec8}

\begin{defs}\label{pr1}
Assume that the plumbing graph $\Gamma$
is a  negative definite connected tree.
We say that $\Gamma$ is {\em almost-rational}
(in short {\em AR})  if there exists
a vertex $j_0\in\calj$ of $\Gamma$ such that replacing $e_{j_0}$ by  some  $e_{j_0}'\leq e_{j_0}$
 we get a rational graph $\Gamma'$.
 In general, the choice of $j_0$ is not unique.
Once  the distinguished vertex $j_0$ is  fixed, we write
$\calj=\{0\}\cup \calj^*$ such that the index $0\in\calj  $
corresponds to this vertex.\end{defs}

\begin{h}\textbf{Examples}\qua
\label{pr2}
The set of {\em AR} graphs is (surprisingly) large.

(1)\qua All rational graphs are {\em AR}. Indeed, take $e_{j_0}'=e_{j_0}$
for any vertex $j_0$.

(2)\qua Any weakly elliptic graph is {\em AR}.  This can be proved as
follows. Take a computation sequence  $x_0,\ldots,x_t$ for
$x_\text{min}$ as in Remark~\ref{artin}.  Since $\Gamma$ is weakly elliptic if
and only if $\chi(x_\text{min})=0$ (see Remark~\ref{artin}, for example),
we get from a similar  discussion like in Remark~\ref{artin} that
$(x_l,b_{j(l)})=1$ all the time   in the computation sequence,
excepting exactly one $l=l_0$  when this  intersection $=2$ (and
$b_{j(l_0)}\in |x_{l_0}|$ as well). Then take $j(l_0)$ for the
distinguished vertex.  If we replace $e_{j(l_0)}$ by a strict
smaller integer, then for some $t'\leq t$, the sequence
$x_0,\ldots, x_{t'}$ will be a computation sequence for the new
fundamental cycle with $(x_l,b_{j(l)})=1$ for any $l<t'$. Hence
$\chi$ of the new fundamental cycle will be 1, a fact which
characterizes the rational singularities.

(3)\qua Any  star-shaped graph  is {\em AR}.  Indeed, first blow down
all the $(-1)$--vertices different from the central vertex (this
transformation preserves the {\em AR} graphs); let this new graph
be $\overline{\Gamma}$. Then take for  $j_0$ the central vertex of
$\overline{\Gamma}$, and take for $-e_{j_0}'$ an integer larger than
the number of adjacent vertices of the central vertex of
$\overline{\Gamma}$. Then the modified graph will become rational. This
follows from the fact that a graph is rational provided that
$-e_j\geq \delta_j$ is satisfied for all its vertices $j$ (see
Section~\ref{4.4}, or \cite{Five}). (In other words, if one takes for the
distinguished vertex of $\Gamma$ the central vertex, and one takes
$-e_{j_0}$ sufficiently large, than one gets a rational graph.)

The same argument shows that
the graphs considered  by Ozsv\'ath and Szab\'o in \cite{OSzP} (and in
Theorem \ref{OSZPI}) are {\em AR}.

(4)\qua The class of {\em AR} graphs is closed while taking subgraphs
and decreasing the Euler numbers $e_j$ (since  the class of rational
graphs is so).

(5)\qua Not every graph is {\em AR}. For example, if $\Gamma$
has two (or more) vertices, both with decoration $-e_j\leq \delta_j-2$,
then $\Gamma$ is not {\em AR}. For example:
\end{h}

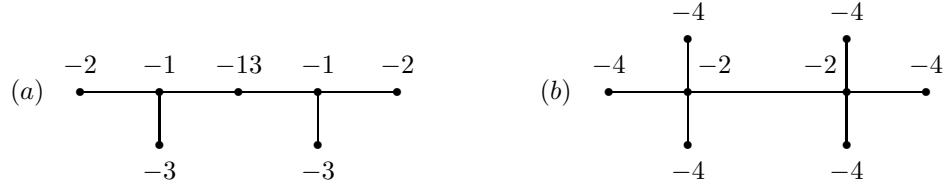
\begin{picture}(375,70)(-20,0)\small
\linethickness{.5pt}
\put(10,40){\circle*{3}}
\put(40,40){\circle*{3}}
\put(70,40){\circle*{3}}
\put(100,40){\circle*{3}}
\put(130,40){\circle*{3}}
\put(40,20){\circle*{3}}
\put(100,20){\circle*{3}}
\put(10,40){\line(1,0){120}}
\put(40,40){\line(0,-1){20}}
\put(100,40){\line(0,-1){20}}
\put(10,50){\makebox(0,0){$-2$}}
\put(40,50){\makebox(0,0){$-1$}}
\put(70,50){\makebox(0,0){$-13$}}
\put(100,50){\makebox(0,0){$-1$}}
\put(130,50){\makebox(0,0){$-2$}}
\put(40,10){\makebox(0,0){$-3$}}
\put(100,10){\makebox(0,0){$-3$}}
\put(-10,40){\makebox(0,0){$(a)$}}

\put(210,40){\circle*{3}}
\put(240,40){\circle*{3}}
\put(300,40){\circle*{3}}
\put(330,40){\circle*{3}}
\put(240,20){\circle*{3}}
\put(240,60){\circle*{3}}
\put(300,20){\circle*{3}}
\put(300,60){\circle*{3}}
\put(210,40){\line(1,0){120}}
\put(240,60){\line(0,-1){40}}
\put(300,60){\line(0,-1){40}}
\put(210,50){\makebox(0,0){$-4$}}
\put(250,50){\makebox(0,0){$-2$}}
\put(290,50){\makebox(0,0){$-2$}}
\put(330,50){\makebox(0,0){$-4$}}
\put(240,10){\makebox(0,0){$-4$}}
\put(300,10){\makebox(0,0){$-4$}}
\put(240,70){\makebox(0,0){$-4$}}
\put(300,70){\makebox(0,0){$-4$}}
\put(190,40){\makebox(0,0){$(b)$}}
\end{picture}

Example (a) was already considered in \cite[Section~4]{OSzP}
(see also Remark~\ref{push}(b)).
Example (b) has the following property: if we delete one of the two
 $(-2)$--vertices, then all the components
of the remaining graph are rational. Still, the graph itself is not {\em AR}.

\begin{thms}\label{isohh}
Assume that the plumbing graph $\Gamma$ is  AR,
 and let $M=M(\Gamma)$ be the oriented plumbed 3--manifold associated with it.
Then for any $[k]\in \Spin^c(M)$,  $HF^+_\text{odd}(-M,[k])=0$, and one has an isomorphism:
$$HF^+_\text{even}(-M,[k])=\bH^+(\Gamma,[k]).$$
In particular,
$$d(M,[k])=\max_{k'\in[k]}\frac{(k')^2+s}{4}=\frac{k_r^2+s}{4}-2\min\chi_{k_r}
=\frac{K^2+s}{4}-2\min_{x\in L}\chi(l'_{[k]}+x);$$
and
$$\chi(HF^+(-M,[k]))=\rank_\Z\bH^+_\text{red}(\Gamma,[k])=\rank_\Z\bH_\text{red}(R_{[k]},\chi_{[k]}).$$

If $\Gamma$ is rational or weakly elliptic then
$$d(M,[k])=\frac{k_r^2+s}{4}=\frac{K^2+s}{4}-2\chi(l'_{[k]}).$$
\end{thms}

\begin{proof}
The proof is similar to the proof of \cite[Theorem~2.1]{OSzP} (see
Theorem~\ref{OSZPI} here),
but the  starting point of the inductive procedure
is the characterization of rational singularities in Theorem~\ref{4.5}.

For any graph $\Gamma$ and vertex $j\in\calj$, we consider the following
two modified graphs:
\begin{itemize}
\item $\Gamma\setminus j $ is obtained from $\Gamma$ by deleting the vertex
$j$ and all the edges adjacent to $j$.

\item $\Gamma^-_j$ is obtained from $\Gamma$ by modifying the decoration
$e_j$ into $e_j-1$.
\end{itemize}
Notice that, in general, $\Gamma\setminus j$ is not connected, write
$\{(\Gamma\setminus j)_c\}_c$
for its connected components. In this case,
$M(\Gamma\setminus j)=\#_cM((\Gamma\setminus j)_c)$,
and all the invariants $HF^+$ and $\bH^+$ are defined in this case as well.
Moreover, one can verify that
 if $\Gamma$ is negative definite, then $\Gamma^-_j$ and all the
components of $\Gamma\setminus j$ are also negative definite.

The proof is based on the following commutative diagram (see
\cite[Section~2]{OSzP}).
We will write $HF^+_\text{o/e}$ for $HF^+_\text{odd/even}$.

{\scriptsize

\begin{picture}(375,70)(0,0)
\linethickness{.5pt}
\put(0,60){\makebox(0,0)[l]{$
HF^+_o(-M(\Gamma\setminus j))   \to
HF^+_e(-M(\Gamma))   \stackrel{A^+}{\to}
HF^+_e(-M(\Gamma^-_j ))   \stackrel{B^+}{\to}
HF^+_e(-M(\Gamma\setminus j))   \stackrel{C^+}{\to}
HF^+_o(-M(\Gamma)) $}}
\put(130,50){\vector(0,-1){30}}
\put(215,50){\vector(0,-1){30}}
\put(300,50){\vector(0,-1){30}}
\put(140,35){\makebox(0,0){$T^+_{\Gamma}$}}
\put(227,35){\makebox(0,0){$T^+_{\Gamma^-_j}$}}
\put(313,35){\makebox(0,0){$T^+_{\Gamma\setminus j}$}}
\put(115,10){\makebox(0,0)[l]{$
\bH^+(\Gamma) \ \ \ \ \   \stackrel{\bA^+}{\longrightarrow} \ \ \  \
\bH^+(\Gamma^-_j)   \ \ \ \ \ \stackrel{\bB^+}{\longrightarrow} \ \ \ \
\bH^+(\Gamma\setminus j)$}}
\end{picture}

}

One has the following properties:

(i)\qua the first line is exact (see \cite[Section~2]{OSzP} and 
\cite[Section~10.12]{OSz7});

(ii)\qua $\bB^+\circ \bA^+=0$, and $\bA^+$ is injective (see 
\cite[Section~2.10]{OSzP}).

\textbf{Step 1}\qua  First, let us prove that  for any rational graph $\Gamma$ one has:

(a)\qua the morphism $T^+_\Gamma\co HF^+_e(-M(\Gamma))\to \bH^+(\Gamma)$ is an
isomorphism,

(b)\qua $HF^+_o(-M)=0$, in particular,
$HF^+_\text{red}(-M)=0$
by Proposition~\ref{hbhb} and Theorem~\ref{4.5}.

If $s$ denotes the number of vertices of $\Gamma$, then this property is true
for $s\leq 3$. Indeed, these cases are covered by \cite{OSzP} (see also
Theorem~\ref{OSZPI}).
By induction, assume that (a) and (b) are true for any rational graph
with number of vertices $<s$. Then, by the general theory (see
\cite[Section~7.2]{OSz}), it is true even for those non-connected
graphs whose connected components are rational with total number of
vertices  $<s$.

Consider now an arbitrary rational graph $\Gamma$ with $s$ vertices. It is
well-known that if $G$ is a rational graph and $j$ is one of its vertices,
, then $G^-_j$  is also rational. Applying this fact several times
starting from $\Gamma$, we get a graph which has the property $-e_j\geq
\delta_j$ for all of its vertices. For such a graph we already know that
(a) and (b) are true (see Theorem~\ref{OSZPI}).  Hence, going back
along this sequence of graphs, we are in the following situation at each
inductive step:

(iii)\qua $T^+_{\Gamma^-_j}$ and $T^+_{\Gamma\setminus j}$ are isomorphisms;

(iv)\qua $HF^+_o(-M(\Gamma^-_j))=HF^+_o(-M(\Gamma\setminus j))=0$.

Then the following fact follows easily from the above diagram and
(i)--(iii):

(v)\qua $\ker \bB^+=\im \, \bA^+$.

Moreover,

(vi)\qua $B^+$ is onto.

This follows by a similar argument as the proof of 
\cite[Proposition~2.8]{OSzP}:
Since both $\Gamma$ and $\Gamma^-_j$ are negative definite,
the two cobordisms inducing $A^+$ and $B^+$ are both negative
definite, hence the third cobordism (inducing $C^+$) is not.
In particular, it induces the trivial map on $HF^\infty$.
But $HF^+_\text{red}(-M(\Gamma\setminus j))=0$ (by induction, since $\Gamma\setminus j$ has
$<s$ vertices), it follows that $C^+=0$.

This induction ends the proof of step 1.

\textbf{Step 2}\qua Assume that $\Gamma$ is an {\em AR} graph
with distinguished vertex $j_0$.  Then all the components of
$\Gamma\setminus j_0$ are rational, hence  $T^+_{\Gamma\setminus
j_0}$ is an isomorphism. Moreover, by the same argument as in
(vi), we get that $C^+=0$ again. Then one can use a similar
inductive procedure as above connecting the graph $\Gamma$ to a
rational graph  by a sequence of {\em AR} graphs, at each time
decreasing $e_{j_0}$ by one.

For the last formulas, see Proposition~\ref{hbhb}, Equation~\ref{3.9}(4),
Section~\ref{luj2} and Remark~\ref{re}(c).
\end{proof}

\begin{h}\textbf{Remarks}\qua
\label{push} The above proof can be pushed further
to obtain the following:

(a)\qua By the same inductive  proof, one can obtain the same result
as in Theorem~\ref{isohh} for a slightly larger class than {\em AR}
graphs. Namely, assume that $\Gamma$ has the following property:
there exists a vertex $j_0$ such that (i) all the components of
$\Gamma\setminus j_0$ are rational; and (ii) replacing $e_{j_0}$
by some $e_{j_0}'\leq e_{j_0}$ we get an {\em AR} graph. Then the
statements of Theorem~\ref{isohh} are valid. Notice that the graph
\ref{pr2}(5)(b) satisfies the above requirement (but it is not {\em
AR}). (Nevertheless, see also Remark~\ref{pr6}(c).)

(b)\qua  Assume that  $\Gamma$ has the following property: there
exists a vertex $j_0$ such that replacing $e_{j_0}$ by some
$e_{j_0}'\leq e_{j_0}$ we get an {\em AR} graph. Then
$\Gamma\setminus j_0$ is an {\em AR} graph, and the proof of
Theorem~\ref{isohh} is still valid excepting (in step 2) the triviality of
$C^+$. Hence, for such graphs we get an isomorphism
$T^+\co HF^+_\text{even}(-M,[k])\to \bH^+(\Gamma,[k])$, but the vanishing
of $HF^+_\text{odd}(-M,[k])$ does not follow. (This result and argument
is the perfect analog of \cite[Theorem~2.2]{OSzP}.)

Notice that the graph \ref{pr2}(5)(a) satisfies the above
requirement. Its homology $HF^+$ is computed in \cite[Remark~4.3]{OSzP},
and $HF^+_\text{odd}$ is  non-zero indeed.\end{h}

\section{{\em AR} graphs and the computation of
$\bH(R_{k_r},\chi_{k_r})$}\label{sec9}

In this section we assume that $\Gamma$ is an {\em AR} graph with
a distinguished vertex $0\in\calj$. We will use the notation
$\calj^*=\calj\setminus \{0\}$. $[k]$ will denote  a fixed
$\spin^c$--structure with its distinguished
 characteristic element $k_r=K+2l'_{[k]}$.
Similarly as in Section~\ref{sec8}, in order to simplify the notation,
sometimes  we write $k$ for $k_r$. For each   $i\geq 0$,  we
consider the cycles $x(i)$, defined in  Lemma~\ref{genlem3}, associated
with the distinguished vertex $0\in \calj$ and $\spin^c$--structure
$[k]$.

\begin{lems}\label{pr3} Assume that $\Gamma$ is a AR graph and fix $i\geq 0$.
Then the following
facts hold.

{\rm(a)}\qua  Consider an arbitrary
$x\in L$ with $m_0(x)=i$. Then $\chi_k(x)\geq \chi_k(x(i))$.

{\rm(b)}\qua Consider any $y>x(i)$ with $m_0(y)=m_0(x(i))=i$.   Then there exists at least one
$b_j\in|y-x(i)| $, for $j\in\calj^*$, such that $\chi_k(y-b_j)\leq \chi_k(y)$.

{\rm(c)}\qua Consider a  computation sequence which connects $x_0=x(i)$ and
$x_t=x(i+1)$, see Lemma~\ref{genlem4}. Then $(x_l+l'_{[k]},b_{j(l)})=1$,
or equivalently, $\chi_k(x_{l+1})=\chi(x_l)$ for any  $0<l<t$. In
particular, $\chi_k(x(i)+b_0)=\chi_k(x(i+1))$.
\end{lems}

\begin{proof} (a)\qua
Write $x=x(i)-y_1+y_2$ with $y_r\geq 0$,  $b_0\not\in |y_r|$ for $r=1,2$, and
$(*)$ \ $|y_1|\cap |y_2|=\emptyset$. Then
$\chi_k(x)=\chi_k(x(i)-y_1)+\chi(y_2)+(y_1,y_2)-(x(i)+l'_{[k]},y_2)
\geq \chi_k(x(i)-y_1)$.  Indeed,
 $(y_1,y_2)\geq 0$ because of  $(*)$,  and  $-(x(i)+l'_{[k]},y_2)\geq 0$
from the definition of $x(i)$.
  Finally, $\chi(y_2)\geq 0$, since  the subgraph supported by
$|y_2|$  can be consider as the subgraph of the modified graph $\Gamma'$,
which is  rational,
hence the subgraph supported by $|y_2|$ itself is rational.
On the other hand,  by Lemma~\ref{genlem5},   $\chi(x(i)-y_1)\geq \chi(x(i))$.

(b)\qua Assume that for any $b_j\in |y-x(i)|$  one has
$\chi_k(y-b_j)>\chi_k(y)$.
This is equivalent with $(b_j,y+l'_{[k]})\geq 2+b_j^2$. Since $(b_j,x(i)+l'_{[k]})\leq 0$,
we get that $(b_j,y-x(i))\geq 2+b_j^2$, or $(b_j,y-x(i)+K)\geq 0$.
Hence $(*)$ $\chi(y-x(i))\leq 0$.   By similar argument as in (a), $y-x(i)$ is supported by a
rational graph, which contradicts $(*)$.

(c)\qua For any $0\leq l\leq t$ define $z_l:=x_l-x(i)$. Then clearly $z_0=0$, $z_1=b_{j(0)}=b_0$,
and $z_t=x(i+1)-x(i)$.
Consider
$(z_l,b_{j(l)})=(x_l+l'_{[k]},b_{j(l)})-(x(i)+l'_{[k]},b_{j(l)})$. Since
$(x_l+l'_{[k]},b_{j(l)})>0$ from the construction of $\{x_l\}$,
$-(x(i)+l'_{[k]},b_{j(l)})\geq 0$ from the definition of $x(i)$, we get
$(z_l,b_{j(l)})>0$.
Notice that $j(l)\in\calj^*$ (for $l>0$)  and $b_0\not\in |z_l|$, hence the values
$(z_l,b_{j(l)})$ will stay unmodified even if we replace our graph $\Gamma$ with the
rational graph $\Gamma'$ (ie if we decrease $e_{v_0}$).
Since $(z_l,b_{j(l)})>0$, we get that (in $\Gamma'$)
the sequence $\{x_l\}_{0<l\leq t}$ can be completed
to a computation sequence connecting $b_{j(0)}$ and $x_\text{min}$, associated with the
canonical characteristic element $K$.
  But then Remark~\ref{artin}  shows that
$(z_l,b_{j(l)})=1$ for any $0<l<t$.
Clearly, this is true in $\Gamma$ as well.
 Using again the above identity we get that
$(x_l+l'_{[k]},b_{j(l)})=1$  for any $0<l<t$, or equivalently
$\chi_k(x_{l+1})=\chi_k(x_l)$ (see the proof of Lemma~\ref{genlem2}).
\end{proof}

Our goal is the identification of the graded tree
$(R_{k},\chi_{k})$ (where $k=k_r$). First we identify the local minimum
points $\calv_1=\{v\in\calv(R_k):\, \delta_v=1\}$ of $\chi_{k}$,
see Section~\ref{2.2}. Fix one $v\in \calv_1$,  and set $\chi_{k}(v)=n_v$.
Let $C_v$ be the connected component of $\overline{L}_{k,\leq n_v}$
corresponding to $v$.

\begin{lems}\label{pr4}
{\rm(a)}\qua $C_v$ has a unique maximal element  $x_v$ (with respect to
$\leq$);

{\rm(b)}\qua  $x_v \in S_{[k]}$;

{\rm(c)}\qua  Set $i_v:=m_0(x_v)$. Then
$ x(i_v)\in  C_v$.
\end{lems}

\begin{proof} Since $v$ is a local minimum point,  $(*)$ $\chi_k(x)=n_v$ for any $x\in C_v$.
Assume that $x,\, x+b_i$ and $x+b_j$ are in $C_v$ (with $i\not=j$) .  By $(*)$ one gets
that $(x+l'_{[k]},b_i)=(x+l'_{[k]},b_j)=1$.  Then by a computation,  $\chi_k(x+b_i+b_j)=\chi_k(x)
-(b_i,b_j)\leq \chi_k(x)$. By the fact that $v$ is a local minimum, one gets
$\chi_k(x+b_i+b_j)=n_v$. In particular,  $x+b_i+b_j\in C_v$ as well.  This is enough for (a).
Indeed, in the presence of two maximal elements, one can connect them by a path in such a way
that $\chi_k$ is constant along the path. Then any subsequence of type $x+b_i,\, x, \, x+b_j$
can be replaced by $x+b_i,\,  x+b_i+b_j,\, x+b_j$. The repeated application
of this  leads to a contradiction.
For (b) write
$\chi_k(x_v+b_j)>\chi_k(x_v)$ for any $j$,  which gives
$(x_v+l'_{[k]},b_j)\leq 0$,
and the minimal property of $x(i_v)$ gives $x(i_v)\leq x_v$.
Notice also that $m_0(x(i_v))=m_0(x_v)$.   Using inductively part (b) of
Lemma~\ref{pr3}, we can construct
a  computation sequence connecting $x(i_v)$ with $x_v$ along which $\chi_k$ is
 non-decreasing.  But since $v$ is a local minimum, $\chi_k$  along this
computation sequence should be constant.  Hence all the cycles of this sequence are
in $C_v$; in particular, $x(i_v)\in C_v$.
\end{proof}

In the next theorem we recover the graded tree $(R_{k},\chi_k)$ from the
$\chi_k$--values $\{\chi_k(x(i))\}_{i\geq 0}$ of the cycles $\{x(i)\}_{i\geq 0}$.
(Recall that $k=k_r$.)

\begin{thms}\label{pr5}
{\rm(a)}\qua  There exists an integer $l$ such that $\chi_k(x(i+1))\geq
\chi_k(x(i))$ for any $i\geq l$.  In particular, by taking $\tau(i):=\chi_k(x(i))$ for any $i\geq 0$,
the graded tree $(R_\tau,\chi_\tau)$ is well-defined (see Example~\ref{2.3}(3)).

{\rm(b)}\qua $(R_{k},\chi_{k})=(R_\tau,\chi_\tau)$.
\end{thms}

\begin{proof}
First we define a map $\psi\co \calv_1(R_k,\chi_k)\to \{x(i)\}_{i\geq
0}$ by $v\mapsto x(i_v)$, where $i_v$ is determined in Lemma~\ref{pr4}.
We claim that $\psi$ is injective.  Indeed, assume that  for $v_1\not=v_2$
one has $i_{v_1}=i_{v_2}$. Part (c) of Lemma~\ref{pr4} guarantees that
$\chi_{k}(v_1)=\chi_k(x(i_{v_1}))=\chi_k(x(i_{v_2}))=\chi_{k}(v_2)$. Let
this integer be $n$.  Consider in $\overline{L}_{k,\leq n}$ the
components $C_{v_r}$  corresponding to $v_r$, $r=1,2$. By part (c) of
Lemma~\ref{pr4},
$x(i_{v_r})\in C_{v_r}$, a fact which contradicts with $C_{v_1}\cap
C_{v_2}=\emptyset$.

Assume now that there exists a sequence of integers $i_1,
i_1+1,\ldots, i_2$ such that $\chi_k(x(i_1))>\chi_k(x(i_1+1))=\cdots
=\chi_k(x(i_2-1))<\chi_k(x(i_2))$.   Set $n=\chi_k(x(i_1+1))$.  Then in
$\overline{L}_{k,\leq n}$ there is a connected  component  $C_v$ which
contains $x(i_1+1)$; let $v$ be the corresponding vertex of $R_{k}$.
By part (a) of Lemma~\ref{pr3}, $v$ is a local minimum point of $R_{k}$.
Since $R_{k}$ has only  a finitely many local minimum points, there
exists $l$ such that for $i\geq l$ such a sequence $i_1,\ldots, i_2$
cannot occur. This, together with $\lim_{i\to \infty}\chi_k(x(i))=\infty$
proves part (a).

Notice also that if $v\in \calv_1(R_k,\chi_k)$, then
$\chi(v)=x(i_v)$ defines a local minimum vertex of
$(R_\tau,\chi_{\tau})$ This follows from the above discussion
and from the properties of the computation sequences connecting the different
$x(i)$'s.

Let $\psi_*\co \calv_1(R_{k})\to \calv_1(R_\tau)$  be induced by $\psi$.
In fact, the previous discussion also show that $\psi_*$  is a bijection  compatible with the gradings.
In order to finish the proof, consider two different  local minimum points
 $u,v\in \calv_1$  of $R_{k}$.  We have to verify that $\chi_{k}(\sup(u,v))=
\chi_\tau(\sup(\psi_*(u),\psi_*(v)))$.  Using the above notation, consider $i_u$ and $i_v$ and assume
$i_u<i_v$.  Then consider one index $i$ with $i_u<i<i_v$ such that
$\chi_k(x(i))\geq \chi_k(x(j))$ for any $i_u\leq j\leq i_v$. Then, in $R_\tau$, $n_{uv}=
\chi_k(\sup(\psi_*(u),\psi_*(v)))$ is exactly $\chi_k(x(i))$.   We show that
the same is true for $\chi_{k}(\sup(u,v))$ in $R_k$.
Indeed, any path connecting $x(i_u)$ and $x(i_v)$ will
have at least one element whose $b_0$--coefficient
is exactly $i$ (since $m_0(u)=i_u<i<i_v=m_0(v)$).   Hence by
Lemma~\ref{pr3} part (a),
$\chi_{k}(\sup(u,v))\geq \chi_k(x(i))$. But the iterated application of
Lemma~\ref{genlem4},
and Lemma~\ref{pr3} part (c), provides a computation sequence connecting
$x(i_u)$ and $x(i_v)$  along which $\chi_k$ is $\leq \chi_k(x(i))$; hence the equality follows.
\end{proof}

\begin{h}\textbf{Remarks}\qua
\label{pr6}
(a)\qua Theorem~\ref{pr5}  shows that for almost rational graphs,
any graded tree $(R_{k},\chi_{k})$ is completely determined by the
values of $\chi_k$ along a very natural infinite computation sequence
(depending on $k$) which contains the elements $\{x(i)\}_{i\geq 0}$.
This sequence can be constructed by iterating Lemma~\ref{genlem4}, ie by
gluing the corresponding sequences  which connect $x(i)$ with $x(i+1)$
for all $i\geq 0$.

In particular, all the important vertices of $R_{k}$ can be represented
by  some special elements  in $L$  which can be arranged in an increasing
linear order (with respect to $\leq$).  (This property  is  extremely
useful in the study  of the weakly elliptic singularities \cite{Ninv},
but here we realize that it is a more general phenomenon.)

(b)\qua The set $\{x(i)\}_{i}$ sometimes is not very economical:  only
some of the $i$'s carry substantial information which will survive
in $(R_\tau,\chi_\tau)$.  For example, for rational singularities,
$\chi(x(i+1))\geq \chi(x(i))$, hence only the information $\chi(x(0))=0$
is preserved in $R_\tau$.

(c)\qua Part (c) of Lemma \ref{pr3} cannot be extended for the class of
graphs considered in Remark~\ref{push}(a). In the case of that  class, one can
still identify $(R_k,\chi_k)$ with some $(R_\tau,\chi_\tau)$, but in
the definition of $\tau$ it is not enough to take only the collection
of cycles $x(i)$ along the infinite computation sequence, but one has
to add some other special cycles as well.

(d)\qua As we already mentioned, one of our interests is the numerical
invariant
$$-\ssw^{OSz}_{M,[k]}=\chi(HF^+(-M,[k]))+\frac{d(M,[k])}{2}.$$
In the case of {\em AR} graphs, this equals
$$\rank_\Z\bH_\text{red}(R_{k_r},\chi_{k_r})-\min\chi_{k_r}+\frac{k_r^2+s}{8}.$$
The term $k_r^2+s$  (from $K^2+s$ and $l'_{[k]}$)  is clear. The
next fact  provides the other term.
\end{h}

\begin{cors}\label{pr7}
For any $\spin^c$--structure $[k]$ consider the
representative $k_r=K+2l'_{[k]}$.  Then
$$\rank_\Z \bH_\text{red}(R_{k_r},\chi_{k_r})-\min\chi_{k_r} =\sum_{i\geq 0}\max
\{ 0, -1+(b_0,x(i)+l'_{[k]})\}.$$
(In fact,
$\rank_\Z\bH_\text{red}(R_{k'},\chi_{k'})$ is independent of the choice
of $k'\in[k]$.)
\end{cors}

\proof
We apply Corollary~\ref{2.10}.  Indeed, using
Lemma~\ref{pr3} part (c),
$$\chi_k(x(i))-\chi_k(x(i+1))=
\chi_k(x(i))-\chi_k(x(i)+b_0)=-1+(b_0,x(i)+l'_{[k]}).\eqno{\qed}$$

\begin{cors}\label{pr7b}
Consider a  normal surface
singularity $(X,0)$, and $\widetilde{X}$ be a fixed resolution of it.
Consider any holomorphic line bundle ${\mathcal L}\in Pic(\widetilde{X})$
on $\widetilde{X}$ with $c_1({\mathcal L})=-l'_{[k]}$ (here $L'\approx
H^2(\widetilde{X},\Z)$ by Poincar\'e duality).
Assume that the resolution graph is AR.  Then

{\rm(a)}\qua $h^1(\widetilde{X},{\mathcal L})
\leq \rank_\Z \bH_\text{red}(R_{k_r},\chi_{k_r})-\min\chi_{k_r}$.

{\rm(b)}\qua If $(X,0)$ is rational then both sides of this inequality
{\rm(a)} are zero.

{\rm(c)}\qua If $(X,0)$ is weakly elliptic, and $[k]=[K]$, then
 $\rank_\Z \bH_\text{red}(R_\text{can})=l$  and $\min\chi_\text{can}=0$ (see
Proposition~\ref{4.7}); in particular $p_g\leq l$.
\end{cors}

This answers positively the inequality part of  the main conjecture \cite{SWI},
modified as in Problem~\ref{quid}(IV) for {\em AR} resolution graphs.
For Gorenstein weakly elliptic
singularities, the equality  $p_g=l$ follows from \cite{Ninv}.
For a generalization of the conjecture \cite{SWI} and the construction of
those line bundles ${\mathcal L}$ for which one expects equality above,
see \cite{NL}.

\begin{proof} We need only  to prove part (a). For this
use the above infinite  computational sequence $\{x_l\}_l$  (see
Remark~\ref{pr6}(a)). First notice that (using the theorem on formal
functions, and the identification of the elements
$x_l=\sum_jn_jb_j\in L$ with cycles $\sum_jn_jE_j$, see
Example~\ref{4.2}), one has
$$H^1(\widetilde{X},{\mathcal L})={\lim_{\leftarrow}} \,
H^1(x_l,{\mathcal L}|_{x_l}).$$ Next, consider the cohomology long
exact sequences associated with the short exact sequences obtained
from
$$0\to {\mathcal O}_{b_{i(l)}}(-x_l)\to {\mathcal
O}_{x_{l+1}}\to {\mathcal O}_{x_{l}}\to  0$$
by tensoring with
${\mathcal L}$. This provides the exact sequence
$$\to H^1({\mathbb P}^1, {\mathcal O}_{{\mathbb
P}^1}(-(b_{i(l)}, x_l+l'_{[k]})))\to H^1({\mathcal
L}|_{x_{l+1}})\to H^1({\mathcal L}|_{x_l})\to 0.$$ Then use
Lemma~\ref{pr3} part (c) and Corollary~\ref{pr7}.
\end{proof}

\section{Example: Lens spaces}
\label{sec10}

\begin{h}\textbf{Notation}\qua We invite the reader to refresh  the
notation of Section~\ref{3.1}. In particular, we recall that
$\{g_j\}_{j\in\calj}$ denotes the dual base in $L'$, see the last
paragraph of Section~\ref{3.1}.\end{h}

\begin{h}\textbf{The plumbing graph of lens spaces}\qua\label{l2} Consider the continued fraction
$$\frac{p}{q}=[k_1,k_2,\ldots,k_s]:=k_1-\cfrac{1}{k_2-\cfrac{1}{\ddots -\cfrac{1}{k_s}}},$$
where $k_j\geq 2$ for any $j$.  Then

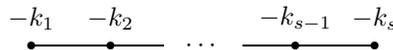
\begin{picture}(200,30)(-50,10)\small
\put(100,20){\circle*{3}}
\put(130,20){\circle*{3}}
\put(200,20){\circle*{3}}
\put(230,20){\circle*{3}}
\put(100,20){\line(1,0){50}}
\put(230,20){\line(-1,0){50}}
\put(165,20){\makebox(0,0){$\cdots$}}
\put(100,30){\makebox(0,0){$-k_1$}}
\put(130,30){\makebox(0,0){$-k_2$}}
\put(230,30){\makebox(0,0){$-k_s$}}
\put(200,30){\makebox(0,0){$-k_{s-1}$}}
\end{picture}

\noindent is the (minimal) plumbing graph $\Gamma$ of the lens space $M=L(p,q)$.
$M$ (and $\Gamma$) can be considered as the link (respectively, the minimal resolution graph)
of the cyclic quotient singularity $\C^2/\Z_p$, where the action is
$\xi*(x,y)=(\xi x,\xi^q y)$. $M$ can also be obtained as a $-p/q$ surgery on the unknot in
$S^3$.

For any $1\leq i\leq j\leq s$ we write the continued fraction
$[k_i,\ldots, k_j]$ as a rational number $n_{ij}/d_{ij}$ with
$n_{ij}>0$ and $\gcd(n_{ij},d_{ij})=1$. We also set $n_{i,i-1}:=1$
and $n_{ij}:=0$ for $j<i-1$. Since
$n_{ij}/d_{ij}=k_i-d_{i+1,j}/n_{i+1,j}$, one gets
$d_{ij}=n_{i+1,j}$ and $n_{ij}=k_in_{i+1,j}-n_{i+2,j}$ for any
$i\leq j$. In fact, any such identity has its symmetric version.
For this it is helpful to notice that $n_{ij}$ is the determinant
of $-B_{(i,j)}$, where $B_{(i,j)}$ is the bilinear form associated
with the graph
$$
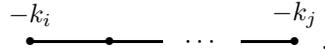
\begin{picture}(120,17)(90,17)\small
\put(100,20){\circle*{3}} \put(130,20){\circle*{3}}
\put(200,20){\circle*{3}} 
\put(100,20){\line(1,0){50}} \put(200,20){\line(-1,0){20}}
\put(165,20){\makebox(0,0){$\cdots$}}
\put(100,30){\makebox(0,0){$-k_i$}}
\put(200,30){\makebox(0,0){$-k_{j}$}}
\end{picture}.$$

This symmetry in the present case reads as:
$n_{ij}=k_jn_{i,j-1}-n_{i,j-2}$ for any $i\leq j$.
 Using these, for any $1\leq j\leq s$ one obtains by induction:
\begin{equation*}
n_{1j}\cdot n_{2s}=p\cdot n_{2j}+n_{j+2,s} \ \ \mbox{and } \ \
n_{js} n_{1,s-1}=p\cdot n_{j,s-1}+n_{1,j-2}.
\tag{1}\end{equation*} Clearly $n_{1s}=p$ and $n_{2s}=q$. We write
$q':=n_{1,s-1}$. Then by (1) $qq'\equiv 1 $ modulo $p$ (and
$0<q'<p$).\end{h}

\begin{h}\textbf{The group $H$, the $\spin^c$ structures, and the elements
$l'_{[k]}$}\qua\label{l3}
Clearly $L'/L=H=\Z_p$, and $[g_s]=g_s+L$ is one of its generators
(with this choice we can use the present formulas in the next
section as well). Indeed,  $[g_j]=[n_{j+1,s}g_s]$ in $H$ $(1\leq
j\leq s)$.

Similarly, the set of $\spin^c$--structures  on $M$ is the set of
orbits $\{[-ag_s]\}_{0\leq a<p}=\{-ag_s+L\}_{0\leq a<p}$ (we
prefer to use this sign, since $-g_s$ is effective). More
precisely, this correspondence is $[k]=K+2(-ag_s+L)$, where $a$
runs from $0$ to $p-1$. In order to emphasize the role of $a$, we
also use the notation  $l'_{[-ag_s]}$ for $l'_{[k]}$.

For any $0\leq a<p$ write
$$l'_{[-ag_s]}=-(a_1g_1+a_2g_2+\cdots+a_sg_s).$$
Since $(l'_{[-ag_s]},b_j)\leq 0$, we get $a_j\geq 0$ for any $j$.
In the next discussion we will clarify the relationship between
the integer $0\leq a<p$ (which codifies $\Spin^c(L(p,q))$) and the
system $E(a):= (a_1,\ldots, a_s)$ (the coefficients of the
corresponding minimal vectors $l'_{[-ag_s]}$).\end{h}

\begin{lemss}\label{l4}
For $i\leq j$ one has
$\sum_{t=i}^j n_{i,t-1}a_t<n_{ij}$ \  and \ $\sum_{t=i}^j
n_{t+1,j}a_t<n_{ij}$.
\end{lemss}

\begin{proof} We verify the second set of inequalities by induction over $j-i$.
If $i=j$, then consider $v:=l'_{[-ag_s]}-b_i\in
l'_{[-ag_s]}+L$. Clearly $(v,b_t)\leq 0$ for any $t\not=i$. Since
$l'_{[-ag_s]} $ is minimal in $(l'_{[-ag_s]} +L)\cap S_\Q$, we get that $(v,b_i)$ should be positive.
Hence $-a_i+k_i>0$.

Next, for the case $j=i+1$, notice that by the inductive step,
$l:=(k_{i}-a_{i})b_{i+1}+b_{i}>0$. Then consider $v:=l'_{[-ag_s]} -l\in l'_{[-ag_s]} +L$.
Then for any $t\not=j$, the value $(v,b_t)\leq 0$ automatically.
Again by the minimality of $l'_{[-ag_s]} $, we get that $(v,b_j)>0$, or $n_{ij}-n_{jj}a_{i}-a_j>0$.
The induction follows similarly.
\end{proof}
In particular, the entries of $(a_1,\ldots,a_s)$ satisfy the system of inequalities:
$$\leqno{\rm(SI)}\quad\left\{ \begin{array}{l}
a_1\geq 0, \ldots , a_s\geq 0\\
n_{i+1,s}a_i+n_{i+2,s}a_{i+1}+\cdots + n_{ss}a_{s-1}+a_s<n_{is}
 \ \mbox{for any $1\leq i\leq s$}.\end{array}\right.$$
By this system one can identify the possible $\spin^c$--structures with the possible
combinations  $(a_1,\ldots,a_s)$.
For this, first notice that the integer $a$ can be recovered from $(a_1,\ldots , a_s)$ by
\begin{equation*}a=
n_{2s}a_1+n_{3s}a_2+\cdots + n_{ss}a_{s-1}+ a_s.\tag{1}\end{equation*}
Indeed, modulo $L$, one has the following  identities:
$$ag_s\equiv -l'_{[-ag_s]}\equiv \sum_{j=1}^s a_jg_j\equiv (\sum_{j=1}^s a_jn_{j+1,s})g_s.$$
Since $\sum a_jn_{j+1,s}<p$ (from (SI), with $i=1$), the above identity follows.

Next, any $0\leq a<p$ determines inductively the entries $a_1,\cdots , a_s$ by the formula
$$a_i=\Big[ \frac{a-\sum_{t=1}^{i-1} n_{t+1,s}a_t}{n_{i+1,s}}\Big] \ \ (1\leq i\leq s).$$
In fact, this also shows that the  set of integral solutions of the system (SI) is
exactly the set of all possible combinations $(a_1,\ldots,a_s)$ associated with the integers
$0\leq a<p$ by the above procedure.

\begin{hh}\label{l5} As a curiosity, we mention that the above system
(SI) can also be interpreted in language of ``generalized''
continued fractions. For any $1\leq i\leq s$ write $r_i:=
n_{is}/n_{i+1,s}$. Then
$$\cfrac{a}{p}=
\cfrac{a_1+\cfrac{ a_2 +\cfrac{  \cdots +\cfrac{a_{s-1}+ \cfrac{a_s}{r_s}}{\cdots}}{r_3}}{r_2}}{r_1}.$$
The inequalities (SI) imply that all the possible fractions in the
above expression are $<1$; and this property guarantees the
uniqueness of the entries $(a_1,\ldots,a_s)$ in this continued
fraction (for any given $0\leq a<p$).\end{hh}

\begin{hh}\label{l6} In fact,
\begin{equation*}
l'_{[-(p-1)g_s]}=g_s+\sum _{i=1}^sb_i= -(\,  (k_1-1)g_1+\sum _{i=2}^{s}(k_i-2)g_i\,).
\tag{1}
\end{equation*}
This follows from the fact that $-(p-1)g_s\equiv g_s$ modulo $L$, and
$g_s+\sum b_j$ is in the unit cube $\sum r_jb_j$, $0<r_j<1$ for any $1\leq j\leq s$.
(Notice that in general it is not true that $l'_{[-ag_s]}$ is in the unit cube.)
(1) can be generalized by the following algorithm and identity
 (the verification is left to the reader).

The set of all possible combinations $E(a)=(a_1,\ldots,a_s)$, $0\leq a<p$, in the order
$E(p-1), \ldots, E(0)$, can be generated inductively as follows.
Start with $E(p-1)=(k_1-1,k_2-2, \ldots, k_s-2)$.
Assume that $E(a)=(a_1,\ldots,a_s)$ is already determined. Then, if
$a_s>0$, then $E(a-1)=(a_1, \ldots,a_{s-1},a_s-1)$. If $a_i=\cdots =a_s=0$, but
$a_{i-1}\not=0$, then
$$E(a-1)=(a_1,\ldots,a_{i-2}, a_{i-1}-1,k_i-1,k_{i+1}-2,\ldots, k_s-2).$$
In particular
$$l'_{[-(a-1)g_s]}-l'_{[-ag_s]}=g_s+\sum_{t=i}^sb_t.$$
\noindent Before we start the list of invariants, we prove the
following identities for $E(a)=(a_1,\ldots, a_s)$. In the sequel
$\{x\}:=x-[x]$ denotes the fractional part of the real number
$x$.\end{hh}

\begin{lems}\label{l7}
$[aq'/p]= \sum_{t=1}^{s}a_tn_{t+1,s-1}$ and $\{aq'/p\}=(\sum_{t=1}^s a_t
n_{1,t-1})/p$.
\end{lems}
\begin{proof} Via  Equation (1) in Lemma~\ref{l4}, we have to show that
$$\sum_{t=1}^{s}a_tn_{t+1,s-1}\leq \frac{n_{1,s-1}}{p} \sum_{t=1}^s a_tn_{t+1,s}<
1+\sum_{t=1}^{s}a_tn_{t+1,s-1}.$$ This follows from the second
part of (1) in Section~\ref{l2} and the first  inequality of
Lemma~\ref{l4} (applied
for $(i,j)=(1,s)$). The second part follows from the expression of
$[aq'/p]$ and the second part of (1) in Section~\ref{l2}.
\end{proof}

\begin{h}\textbf{The Ozsv\'ath--Szab\'o invariant $HF^+(\pm M)$}\qua
\label{l8} Since $\Gamma$ is
rational, by Theorem~\ref{4.5} and Theorem~\ref{isohh}, $HF^+_\text{red}(\pm M)=0$ and
 $HF^+(\pm M,[k])=\calt^+_{\pm d}$, where
$$d:=d(M,[k])=\frac{k_r^2+s}{4}=\frac{K^2+s}{4}-2\chi(l'_{[k]}).$$
Notice that (see eg \cite[Section~7.1]{SWI}):
$$\frac{K^2+s}{4}=\frac{p-1}{2p}-3\cdot \bms(q,p),$$
where $\bms(q,p)$ denotes the Dedekind sum
\[
\bms(q,p)=\sum_{l=0}^{p-1}\Big(\!\Big( \frac{l}{p}\Big)\!\Big)
\Big(\!\Big( \frac{ ql }{p} \Big)\!\Big), \ \mbox{where} \ \
(\!(x)\!)=\left\{
\begin{array}{ccl}
\{x\} -1/2 & {\rm if} & x\in {\R}\setminus {\Z}\\
0 & {\rm if} & x\in {\Z}.
\end{array}
\right.
\]
(In fact, this formula for $K^2+s$ for cyclic quotients goes back
to the work of Hirzebruch.)

Next, we wish to determine $\chi(l'_{[-ag_s]})$.\end{h}

\begin{propss}\label{l9}
For any $0\leq a<p$ one
has
$$\chi(l'_{[-ag_s]}) =\frac{a(1-p)}{2p} +\sum_{j=1}^a
\Big\{\frac{jq'}{p}\Big\}.$$
\end{propss}

\begin{proof} Using the notation of Section~\ref{l6} (in particular
$a_t=0$ for $t\geq i$), we have
$$\chi(l'_{[-(a-1)g_s]})=\chi(l'_{[-ag_s]}
+g_s+\sum_{t=i}^sb_t)=$$ $$\chi( l'_{[-ag_s]})
+\chi(g_s+\sum_{t=i}^sb_t)-(l'_{[-ag_s]}, g_s+\sum_{t=i}^sb_t).$$
Since $ l'_{[-ag_s]}=-\sum_{t=1}^{i-1}a_tg_t$ and $-pg_s=\sum _j
n_{1,j-1}b_j$, and from Lemma~\ref{l7}, one gets
$$-(l'_{[-ag_s]}, g_s+\sum_{t=i}^sb_t) =\sum_ta_t(g_t,g_s)=-\sum_ta_tn_{1,t-1}/p=-\{aq'/p\}.$$
Moreover, $\chi(g_s+\sum_{t=i}^sb_t)=\chi(g_s)+\chi(\sum_{t=i}^sb_t)-1=\chi(g_s)$
by an easy induction.  Notice that
the $b_s$--coefficient of $K$ (ie $(K,g_s)$)
is $-1+(q'+1)/p$ (see \cite[Section~5.2]{SWI}), and $g_s^2=-q'/p$, hence
$\chi(g_s)=(p-1)/(2p)$.  Hence the formula follows by induction over $a$ (since $\chi(l'_{[0]})=0$).
 \end{proof}

\begin{hh} Clearly, with the choice $g_1$ as a generator, by a similar
proof one has
$$\chi(l'_{[-ag_1]}) =\frac{a(1-p)}{2p} +\sum_{j=1}^a
\Big\{\frac{jq}{p}\Big\}.$$\end{hh}

\begin{h}\label{l10} \textbf{The Casson--Walker invariant}\qua  is $\lambda(L(p,q))=p\cdot \bms(q,p)/2$
(see, for example, \cite[Section~5.3]{SWI}).\end{h}

\begin{h}\textbf{The Reidemeister--Turaev torsion}\qua
\label{l11} One can use
for the computation of the sign-refined Reidemeister--Turaev
torsion the formula \cite[Section~5.7]{SWI} (see also
\cite[Section~3.5]{SWI}).
In that article, the canonical $\spin^c$--structure $\sigma_\text{can}$
was defined via
 $c_1(\sigma_\text{can})=-K$, and an arbitrary $\spin^c$ structure
via the action $\sigma=h_{\sigma}\cdot \sigma_\text{can}$ for some
$h_{\sigma}\in H$. This reads as $-c_1(\sigma)=-c_1(\sigma_\text{can})-2h_{\sigma}$,
which in our case corresponds to $k_r=K+2l'_{[k]}$.
In particular, for $[-ag_s]$, $h_{\sigma}$ is the class of $ag_s$ (see also
Section~\ref{4.1}).

Then we apply \cite[Theorem~5.7]{SWI}.  Here
we identify any  character $\chi$ of $H$ by the root of unity $\xi:=\chi([g_s])$. Then we get:
$$\et_{M,[-ag_s]}(1)=\frac{1}{p}\sum\frac{\xi^{-a}}{(\xi-1)(\xi^{q}-1)},$$
where the sum is over all $p$--roots of unity $\xi$ with $\xi\not=1$.
If $a=0$ then
$$\et_{M,[0]}(1)=\frac{1}{p}\sum\frac{1}{(\xi^q-1)(\xi-1)}=\frac{p-1}{4p}-\bms(q,p).$$
For $a>0$, we use the same inductive step as above, namely
$$\et_{M,[-(a-1)g_s]}(1)-
\et_{M,[-ag_s]}(1)=\frac{1}{p}\sum \frac{\xi^{-a}}{\xi^{q}-1}=
\frac{1}{p}\sum \frac{\xi^{-aq'}}{\xi-1}=\Big\{\frac{aq'}{p}\Big\}+\frac{1-p}{2p}.$$
Adding all these identities, and using Proposition~\ref{l9}, one gets
$$\et_{M,[-ag_s]}(1)=\frac{p-1}{4p}-\bms(q,p)-\chi(l'_{[-ag_s]}).$$\end{h}

\begin{h}\textbf{The identity $\ssw^{OSz}_{M,[k]}=\ssw^{TCW}_{M,[k]}$}\qua
\label{l12}
Using the above formulas, one gets
$$\et_{M,[k]}(1)-\frac{\lambda(M)}{|H|}=\frac{d(M,[k])}{2},$$
which proves the  identity $\ssw^{OSz}_{M,[k]}=\ssw^{TCW}_{M,[k]}$
mentioned in Section~\ref{quid}.

Since $\sum_{[k]}\et_{M,[k]}(1)=0$, one also gets
$$-\lambda(M)=\sum_{[k]\in \Spin^c(M)}\ \frac{k_r^2+s}{8},$$
or, equivalently,
$$\sum_{[k]}\chi(l'_{[k]})=\frac{p-1}{4} -p\cdot
\bms(q,p).$$\end{h}

\section{Example:  Seifert 3--manifolds}
\label{sec11}

\begin{h}\label{sei1}  In this section we will use mainly the notation of
\cite{SWII} (which also provides a list of  classical references
for Seifert 3--manifolds).

 Assume that $M$ is a
Seifert rational homology sphere with negative orbifold Euler
number $e$ (equivalently, the corresponding plumbing graph is
negative definite). The minimal plumbing graph $\Gamma$ is
star-shaped with $\nu$ arms; we assume $\nu\geq 3$. (The case
$\nu<3$ is completely covered by the previous section.) $\Gamma$
can be determined from the (normalized) Seifert invariants
$\{(\alpha_l,\omega_l)\}_{l=1} ^\nu$ and $e$ as follows. Here
$0\leq \omega_l<\alpha_l$ for any $1\leq l\leq \nu$, and
 $e=e_0+\sum_l\omega_l/\alpha_l$,
where $e_0$ is the decoration of the central vertex.
 The other decorations can be determined from the continued fraction
$\alpha_l/\omega_l=[k_{1}^\lp, k_{2}^\lp, \cdots, k_{s_l}^\lp]$;
the vertices $\{v_{j}^\lp\}_{j=1}^{s_l}$
situated on the $l^{th}$--arm are decorated by $e_{j}^\lp=-k_{j}^\lp$ in  such a way that the
vertex $v_{1}^\lp$ is  connected with the central vertex.

We denote  by $b_0$ and $b_{j}^\lp \ (1\leq j\leq s_l;\ 1\leq l\leq \nu)$ the
base elements of $L$, respectively by
$g_0$ and $\{g_{j}^\lp\}_{jl}$ the  dual base elements of $L'$ (see
Section~\ref{3.1}).
As we already mentioned (see Example~\ref{pr2}(3)), such a graph $\Gamma$
is {\em AR} with its central vertex distinguished.

We will also use the following notation:
$\ev:=(2-\nu+\sum_l1/\alpha_l)/e$ denotes the
``exponent of $M$'' (in some articles $-\ev $ is the ``log discrepancy of the central vertex'').
For any $l$ we consider $0<\omega_l'<\alpha_l$ with $\omega_l\omega_l'\equiv 1$
modulo $\alpha_l$.

The following identities are well-known (they can be deduced from
\cite[Section~5]{SWI}, for example).
$$(g_0,g_0)=\frac{1}{e};\
(g_0,g^l_{s_l})=\frac{1}{e\alpha_l};\
(g^t_{s_t},g^l_{s_l})=\frac{1}{e\alpha_t\alpha_l} (t\not=l);\
 (g^l_{s_l},g^l_{s_l})=\frac{1}{e\alpha_l^2}-\frac{\omega_l'}{\alpha_l}.$$
$$-(K,g_0)=1+\ev; \ \ -(K,g^l_{s_l})=1+\frac{\ev-\omega'_l}{\alpha_l}.$$
$$\chi(g_0)=\frac{1}{2}+\frac{\ev}{2} -\frac{1}{2e};\ \
\chi(g^l_{s_l})=\frac{1}{2}+\frac{\ev}{2\alpha_l}-\frac{1}{2e\alpha_l^2}.$$\end{h}

\begin{h}\textbf{Notation}\qua
\label{sei2} For any fixed $1\leq l\leq \nu$, we consider similar
notation as in the previous section.
For any $1\leq i\leq j\leq s_l$, we write $n_{ij}^\lp/d_{ij}^\lp:=[k_i^\lp,\ldots, k_j^\lp]$
(with $n_{ij}^\lp>0$ and $\gcd(n_{ij}^\lp,d_{ij}^\lp)=1$); and we set
$n_{i,i-1}^\lp:=1$ and $n_{ij}^\lp:=0$ for $j<i-1$.
Clearly they satisfy the identities  considered  in Section~\ref{l2}.
 Moreover, $\alpha_l=
n_{1,s_l}^\lp$, $\omega_l=n_{2,s_l}^\lp$,  and
$\omega_l'=n_{1,s_l-1}^\lp$.\end{h}

\begin{h}\textbf{The group $H$}\qua
\label{sei3} In the abelian group $H=L'/L$ one has
$$[g_j^\lp]=[n_{j+1,s_l}^\lp g_{s_l}^\lp].$$
In fact,  one has the following (Abelian) presentation (see \cite{Neu}):
\begin{equation*}
H=\langle g_0, g_{s_1}^{1}, g_{s_2}^{2},\ldots g_{s_\nu}^{\nu}\ |\
e_0g_{0}+\sum_{l=1}^\nu \omega_lg_{s_l}^\lp =0,
 \ \ g_{0}=\alpha_lg_{s_l}^\lp\  \mbox{for all $l$} \rangle.
\end{equation*}
 The order of $H$ is $-e\alpha_1\cdots \alpha_\nu$,
and the order $\xo $ of the subgroup $\langle g_{0} \rangle$ is
$\xo=-e\alpha$, where
$\alpha:=\mbox{lcm}(\alpha_1,\ldots,\alpha_\nu)$ (see
\cite{Neu}).\end{h}

\begin{h}\textbf{$\Spin^c$--structures and representatives $l'_{[k]}$}\qua
\label{sei4}
We fix a $\spin^c$--structure $[k]=K+2(l'_{[k]}+L)$, or
equivalently,  the  orbit $l'_{[k]}+L$ in $L'/L$. We assume that
the distinguished representative $l'_{[k]}\in L'$ has the
following form
\begin{equation*}
l'_{[k]}=-a_0g_0-\sum_{j,l}a_{j}^\lp g_j^\lp.
\tag{1}\end{equation*}
For any fixed $1\leq l\leq \nu$, we associate with the collection
$(a_1^\lp,\ldots, a_{s_l}^\lp)$ the integer (see Section~\ref{l3}):
\begin{equation*}
a_l:=\sum_{t\geq 1}n_{t+1,s_l}^\lp a_t^\lp. \tag{2}\end{equation*}
The integers $\{a_0;a_j^\lp\}_{lj}$, respectively
$\{a_0;a_\lp\}_l$ satisfy the following inequalities:\end{h}

\begin{props}\label{sei5} Consider $l'_{[k]}$ and the
notation {\rm(1)} and {\rm(2)} above. Then
$$\leqno{\rm(SI)}\qquad \left\{ \begin{array}{ll}
a_0\geq 0;, a_j^\lp\geq 0 &
(1\leq l\leq \nu, \ 1\leq j\leq s_l);\\
\widetilde{a}^\lp_j:=\sum_{t\geq j} n_{t+1,s_l}^\lp  a_t^\lp <n_{js_l}^\lp &
(1\leq l\leq \nu, \ 1\leq j\leq s_l); \\
1+a_0+ie_0+\sum_l\Big[\frac{i\omega_l+a_l}{\alpha_l}\Big]\leq 0 & \mbox{for any $i>0$}.
\end{array}\right.$$
In particular, the set of integers $(a_0,a_1,\ldots, a_\nu)$ satisfy
$$\leqno{\rm(SI_\text{red})}\qquad \left\{ \begin{array}{ll}
a_0\geq 0;, \ \alpha_l> a_l \geq 0 &
(1\leq l\leq \nu);\\
1+a_0+ie_0+\sum_l\Big[\frac{i\omega_l+a_l}{\alpha_l}\Big]\leq 0 & \mbox{for any $i>0$}.
\end{array}\right.$$
Moreover, there is a one-to-one correspondence between the integer
solutions of {\rm(SI)}, respectively of ${\rm(SI_\text{red})}$,   provided  by
the correspondence {\rm(2)} above (see also Section~\ref{l3}).
\end{props}
\begin{proof}
The first two set of inequalities follow from (SI) in Lemma~\ref{l4}. For
the last one, we proceed as follows. For any $i>0$ we define
$y(i):=ib_0+\sum_{l,j} u_j^l b_j^l \in L,$
where
$$u_1^l:=\Big[\frac{in_{2s_l}^l+\widetilde{a}^l_1}{n_{1s_l}^l}\Big],\
u_2^l:=\Big[\frac{u_1^ln_{3s_l}^l+\widetilde{a}^l_2}{n_{2s_l}^l}\Big],\
\cdots\ ,
u_{s_l}^l:=\Big[\frac{u^l_{s_l-1}+\widetilde{a}^l_{s_l}}{n_{s_ls_l}^l}\Big].$$
Consider $v:=l'_{[k]}-y(i)$. The point is that $(v,b_j^l)\leq 0$ for any
$1\leq l\leq \nu$ and $1\leq j\leq s_l$.  Indeed, if $j<s_l$,
from the definition of $u_j^l$ one has
$u^l_jn^l_{j,s_l} \leq u_{j-1}^ln_{j+1,s_l}^l+\widetilde{a}^l_j,$
or, using Section~\ref{l2} (namely
$n_{j,s_l}^l=k_jn_{j+1,s_l}^l-n_{j+2,s_l}^l$) and
$\widetilde{a}_j^l=\widetilde{a}_{j+1}^l+n_{j+1,s_l}^la_j^l$,
$$-a_j^l-u_{j-1}^l+k_ju_j^l\leq
  \frac{u_j^ln_{j+2,s_l}^l+\widetilde{a}^l_{j+1}}{n^l_{j+1,s_l}},$$
hence $-a_j^l-u_{j-1}^l+k_ju_j^l\leq$ the integral  part of the
right hand side, which equals $u^l_{j+1}$. In other words,
$(v,b^l_j)=-a_j^l-u_{j-1}^l+k_ju_j^l-u^l_{j+1}\leq 0$. (The case
$j=s_l$ is easy.)

Since $v=l'_{[k]}-y(i)\in l'_{[k]}+L$ and $y(i)>0$, by the minimality of $l'_{[k]}$ in
$(l'_{[k]}+L)\cap S_\Q$ we deduce that $v\not\in S_\Q$. Since $(v,b_j^l)\leq 0$ for any $l$
and $j$, we get that $(v,b_0)>0$. This is equivalent to
$-a_0-ie_0-\sum_lu_1^l>0$, hence the third inequality of $(SI)$ follows.

For ${\rm(SI_\text{red})}$, notice that $a_l=\widetilde{a}_1^l<
n_{1s_l}^l=\alpha_l$ from (SI). The correspondence
between the integral solutions of the two system of inequalities is
clarified in Section~\ref{l3};
namely, is given by $E(a_l)=(a^l_1, \ldots, a^l_{s_l})$.
\end{proof}
The point is that the set of solutions of (SI) (or,
equivalently, of ${\rm(SI_\text{red})}$) corresponds exactly to the set of
all possible $\spin^c$--structures. In other words, if
$\{a_0;a_j^l\}_{jl}$ satisfies ${\rm(SI_\text{red})}$, then it is the set of
coefficients of some $l'_{[k]}$. This follows from the next fact.

\begin{props}\label{sei6}
For any $c_0\geq 0$ and $\alpha_l>c_l\geq 0$ consider
$E(c_l)=(c_1^l,\ldots , c^l_{s_l})$ as in Section~\ref{l3}. This defines
$l':= -c_0g_0-\sum_{l,j}c_j^lg_j^l\in S_\Q$. If $l'$ is not  minimal in
$(l'+L)\cap S_\Q$ then for some $i>0$ one has
$$1+c_0+ie_0+\sum_l\Big[\frac{i\omega_l+c_l}{\alpha_l}\Big]>0.$$
\end{props}
\begin{proof}
Define the integers $\widetilde{c}_j^l$, $u_j^l$ and the cycle $y(i) $
(for any $i>0$)
in a similar way as in (SI) and in previous proof, by replacing $a_j^l$ by
$c_j^l$. Then the argument of the previous proof repeated gives that
$y(i)$ satisfies

(a)\qua $m_0(y(i))=i$, and

(b)\qua $(l'-y(i),b_j^l)\leq 0$ for any $l$ and $j$.

The point is
that for any fixed $i>0$, $y(i)$ is the unique maximal effective
cycle which satisfy both (a) and (b).
This follows from the very
definition of the integers $u^l_j$,
its proof is similar to the proof of  Proposition~\ref{sei11}  and it
is left to the reader.

Now, assume that
$l'=l'_{[k]}+x$ for some effective non-zero cycle $x\in L$. Set $i:=m_0(x)$.
Since $x$ satisfy both (a) and (b) for this $i$, from the maximality of $y(i)$
one gets $y(i)\geq x$. Consider now $l'-y(i)=l'_{[k]}+x-y(i)$. Since
$(l'-y(i),b^l_j)\leq 0$ one gets that $(l'_{[k]}+x-y(i),b_j^l)\leq 0$ as well.
Moreover, $(l'_{[k]}+x-y(i),b_0)\leq (l'_{[k]},b_0)\leq 0$  (the first inequality
from $b_0\not\in|y(i)-x|$ and $x-y(i)\leq 0$, the second from $l'_{[k]}\in S_\Q$).
Hence $l'_{[k]}+x-y(i)\in S_\Q$ with $x-y(i)\leq 0$. From the minimality of
$l'_{[k]}$ one gets $x=y(i)$. Hence $(b_0,l'-y(i))=(b_0,l'_{[k]})\leq 0$, or
$-c_0-ie_0-\sum_l[i\omega_l+c_l/\alpha_l]\leq 0$.
\end{proof}
\begin{cors}\label{sei7}
The following sets are in one-to-one correspondence:

{\rm(1)}\qua $\Spin^c(M)$, or equivalently, the set of distinguished representatives $k_r=K+2l'_{[k]}$;

{\rm(2)}\qua the set of integral solutions $\{a_0;a_j^l\}_{lj}$ of {\rm(SI)};

{\rm(3)}\qua  the set of integral solutions $\{a_0;a_1,\ldots,a_\nu\}$ of
${\rm(SI_\text{red})}$.

The set $\{a_0;a_j^l\}_{lj}$ is the coefficient set of $l'_{[k]}$,
$\{a_0,\ldots,a_\nu\}$ is obtained from $\{a_0;a_j^l\}_{lj}$
by Equation~(2) of Section~\ref{sei4}. Notice also that in $L'$ one has
$$-l'_{[k]}\equiv a_0g_0+\sum_{l,j} a_j^lg_j^l\equiv a_0g_0+
\sum_l a_lg^l_{s_l}\ \ (\mbox{modulo}\ L).$$ In particular, the
cardinality of all these sets is the same
$|H|=|e|\alpha_1\cdots\alpha_\nu$. The author knows no direct
proof of this fact for {\rm(2)} or {\rm(3)}.\end{cors}

\begin{h}\textbf{The values $\chi(l'_{[k]})$}\qua
\label{sei8} Assume that $[k]$
 corresponds (by the correspondence in Corollary~\ref{sei7}) to
the set $\{a_0;a_1,\ldots ,a_\nu\}$. Assume that $a_l>0$ for some
fixed $1\leq l\leq \nu$.  Let $[k^-_l]$ be the $\spin^c$--structure
corresponding to
$$\{a_0;a_1,\ldots,a_{l-1},a_l-1,a_{l+1},\ldots, a_\nu\}$$
(clearly, this set of integers also satisfies
${\rm(SI_\text{red})}$). Then, by a similar method as in (the proof of)
Proposition~\ref{l9},  we get
$$\chi(l'_{[k^-_l]})=\chi(l'_{[k]})+\chi(g^l_{s_l})+\frac{a_0}{e\alpha_l}+
\sum_{t\not=l}\frac{a_t}{e\alpha_t\alpha_l} +\sum_{t\geq 1} a_t^l(g^l_t,g^l_{s_l}).$$
By a calculation (using Lemma~\ref{l7} and  Proposition~\ref{l9}):
$$\sum_{t\geq 1} a_t^l(g^l_t,g^l_{s_l})=
\sum_{t \geq 1} a_t ^l
\Big(\frac{n^l_{t+1,s_l}}{e\alpha_l^2}-\frac{n^l_{1,t-1}}{\alpha_l}\Big)=
\frac{a_l}{e\alpha_l^2}-
\Big\{\frac{a_l\omega'_l}{\alpha_l}\Big\}.$$
Therefore,
$$\chi(l'_{[k^-_l]})-\chi(l'_{[k]})=\chi(g^l_{s_l})+\frac{1}{e\alpha_l}\Big( a_0+\sum_{t=1}^\nu
\frac{a_t}{\alpha_t}\Big)-
\Big\{\frac{a_l\omega'_l}{\alpha_l}\Big\}.$$ Hence, by induction
(by decreasing the coefficients $(a_1,\ldots,a_\nu)$ to
$(0,\ldots,0)$), $\chi(l'_{[k]})$ can be computed in terms of
$\chi(l'_{[-a_0g_0]})$. But, using the identities in Section~\ref{sei1}, one
has:
$$-\chi(l'_{[-a_0g_0]})=\frac{1}{2}(\, -a_0(K,g_0)+a_0^2g_0^2\,)=
\frac{a_0^2}{2e}+\frac{a_0}{2}(1+\ev).$$
Summing up these identities, and using $\widetilde{a}:=a_0+\sum_la_l/\alpha_l$, one gets
$$-\chi(l'_{[k]})=\sum_{l=0}^\nu \frac{a_l}{2}+\frac{\ev\widetilde{a}}{2}+\frac{\widetilde{a}^2}{2e}-
\sum_{l=1}^\nu\sum_{i=1}^{a_l}\Big\{\frac{i\omega_l'}{\alpha_l}\Big\}.$$
For another expression see ${\rm(r*)}$ in Section~\ref{sei15}.\end{h}

\begin{h}\textbf{The invariant $k_r^2+s$}\qua
\label{sei9} Recall that
$k_r^2+s=K^2+s-8\chi(l'_{[k]})$.
$\chi(l'_{[k]})$ is computed above in Section~\ref{sei8}, and $K^2+s$ in
\cite[Section~5.5]{SWI}:
$$K^2+s=\ev^2e+e+5-12\sum_{l=1}^\nu \bms
(\omega_l,\alpha_l).$$\end{h}

\begin{h}\textbf{The cycles $x(i)$ for $i\geq 0$}\qua
\label{sei10} We fix a $\spin^c$--structure $[k]$
which corresponds to the set $\{a_0;a_1,\ldots,a_\nu\}$ by the
correspondence in Corollary~\ref{sei7}. Recall that the graded root $(
R_{k_r},\chi_{k_r})$, in particular the $\Z[U]$--module
$\bH(R_{k_r},\chi_{k_r})$, can be completely recovered from the
cycles $x(i)$ ($i\geq 0$) (see Sections \ref{cs} and \ref{sec9}). The next result
determines these cycles in terms of the Seifert invariants of $M$
and the integers $\{a_0,\ldots, a_\nu\}$.\end{h}

\begin{props}\qua
\label{sei11} For any $l$ consider $E(a_l)=(a_1^l,\ldots, a^l_{s_l})$
(see Proposition~\ref{sei5} and Section~\ref{l3}), and define
$\widetilde{a}^l_j:=\sum_{t\geq j} n_{t+1,s_l}^la_t^l$ (for $1\leq j\leq s_l$)
(see $\mathrm{(SI)}$ in Proposition~\ref{sei5}).

{\rm(1)}\qua Define the integers $\{v^l_j\}$ $(1\leq l\leq\nu,\ 1\leq j\leq s_l)$ inductively by
$$v^l_1:=\Big\lceil \frac{i\omega_l-a_l}{\alpha_l}\Big\rceil=
\Big\lceil \frac{in^l_{2,s_l}-\widetilde{a}^l_{1}}{n^l_{1,s_l}}\Big\rceil;\ \
v^l_j:=
\Big\lceil \frac{v^l_{j-1}n^l_{j+1,s_l}-\widetilde{a}^l_j}{n^l_{j,s_l}}\Big\rceil \ (1<j\leq s_l).$$
Set $z(i):= ib_0+\sum_{l,j}v_j^lb_j^l$. Then $(l'_{[k]}+z(i),b^l_j)\leq 0$ for any $l$ and $j$.

{\rm(2)}\qua Assume that an effective cycle
$z(i):= ib_0+\sum_{l,j}v_j^lb_j^l$ (where $\{v^l_j\}_{lj}$ are some integers) satisfies
$(l'_{[k]}+z(i),b^l_j)\leq 0$ for any $l$ and $j$. Then
$$v^l_1\geq \Big\lceil \frac{i\omega_l-a_l}{\alpha_l}\Big\rceil \ \ \mbox{and} \ \
v^l_j\geq
\Big\lceil \frac{v^l_{j-1}n^l_{j+1,s_l}-\widetilde{a}^l_j}{n^l_{j,s_l}}\Big\rceil \ (1<j\leq s_l).$$

{\rm(3)}\qua In particular, $z(i)$ defined in (1) equals $x(i)$,  the
minimal effective cycle which satisfies $(l'_{[k]}+x(i),b^l_j)\leq
0$ for any $l$ and $j$.
\end{props}
\begin{proof} The proof of (1) has the same spirit as the proof of
Proposition~\ref{sei5}.
From the definition of $v^l_j$ one gets $n^l_{j,s_l}v^l_j\geq
v^l_{j-1}n^l_{j+1,s_l}-\widetilde{a}^l_j$,
which is equivalent to
$$k^l_jv^l_j-v^l_{j-1}+a^l_j\geq \frac{ v^l_jn^l_{j+2,s_l}-\widetilde{a}^l_{j+1}}{n^l_{j+1,s_l}}.$$
Now, using the definition of $v^l_{j+1}$  for the right hand side, one gets
$-a^l_j+v^l_{j-1}-k^l_jv^l_j+v^l_{j+1}\leq 0$.
(Here, if $j=1$  then $v^l_{j-1}:=i$.)
For (2), notice that $(l'_{[k]}+x(i),b^l_{s_l})\leq 0$ is equivalent with the wanted
inequality for $j=s_l$.  In general, for arbitrary $j$, the inequality follows from
$(l'_{[k]}+x(i),v_j)\leq 0$, where $v_j=\sum_{t\geq j} n^l_{t+1,s_l}b^l_t$.
\end{proof}

\begin{h}\textbf{The values $\tau(i)$ ($i\geq 0$)}\qua
\label{sei12} Recall that $\tau(i):=\chi_{k_r}(x(i))$ ($i\geq 0
$).
If $i=0$ then $x(0)=0$, hence $\tau(0)=0$.  From Lemma~\ref{pr3} part (c)
and Proposition~\ref{sei11}, one gets that
for $i\geq 0$
$$\tau(i+1)-\tau(i)=1-(l'_{[k]}+x(i),b_0)=1+a_0-ie_0+\sum_l \Big[
\frac{-i\omega_l+a_l}{\alpha_l}\Big].$$ In particular, for $i=0$
one has $\tau(1)-\tau(0)=1+a_0\geq 1$. In general,
$$\tau(i)=\sum_{t=0}^{i-1}\Big(1+a_0-te_0+\sum_l \Big[
\frac{-t\omega_l+a_l}{\alpha_l}\Big]\, \Big).$$\end{h}

\begin{h}\textbf{The Ozsv\'ath--Szab\'o invariant}\qua
\label{sei13} Using the general results
Proposition~\ref{hbhb}, Theorem~\ref{isohh} and Theorem~\ref{pr5}, and
$$\min\chi_{k_r}=\min _{i\geq 0}\tau(i),$$
one gets
$$HF^+_\text{even}(-M,[k])=\bH^+(\Gamma,[k])=\bH(R_\tau,\chi_\tau)[-\frac{k_r^2+s}{4}].$$
Clearly, from its very  definition, the graded root $(R_\tau,\chi_\tau)$ is completely
determined from the values $\{\tau(i)\}_{i\geq 0}$. Moreover,
the numerical invariants can be determined as follows:
$$\pm d(\pm M,[k])=\frac{k_r^2+s}{4}-2\min _{i\geq 0}\tau(i),$$
and (using Corollary~\ref{2.10} as well):
$$\chi(HF^+(-M,[k]))=
\min_{i\geq 0}\tau(i)+\sum_{i\geq 0}\max\Big\{ 0,\,
-1-a_0+ie_0-\sum_{l=1}^\nu \Big[
\frac{-i\omega_l+a_l}{\alpha_l}\Big]\, \Big\}.$$ (Since
$e=e_0+\sum\omega_l/\alpha_l<0$, this sum is finite.)\end{h}

\begin{h}\textbf{Remark}\qua
\label{sei14}  Assume that $[k]=[K]$, ie $l'_{[k]}=0$
(or, equivalently, $a_0=\cdots =a_\nu=0$). Then the above sum becomes
(the Dolgachev--Pinkham) invariant
$$DP:=\sum _{i\geq 0} \max\Big\{0, -1+ie_0
-\sum_{l}\Big[ \frac{-i\omega_l}{\alpha_l}\Big ] \Big\},$$
and the above  identity
$$\chi(HF^+(-M,can))-\min\chi_\text{can} =DP.$$
This is interesting for two reasons.

(1)\qua If $M$  is the link of a normal surface singularity with
geometric genus $p_g$, then Corollary~\ref{pr7b} reads as  $p_g\leq DP$. In
fact, if $M$ is a rational homology sphere, and it is the link of
a singularity which admits a good $\C^*$--action (a fact which
automatically implies that $M$ is a Seifert 3--manifold), then by
results of Dolgachev and Pinkham, $p_g=DP$ (see eg
\cite{Pinkh}). It is remarkable, that in that context, the
expression $DP$ is deduced by a rather different argument (sitting
in algebraic geometry).

This inequality (and identity for  singularities with
good $\C^*$--action) is compatible with the conjecture \cite{SWI},
see also \cite{SWII} for this special situation.

(2)\qua The identity $\ssw^{OSz}_{M,\text{can}}=\ssw^{TCW}_{M,\text{can}}$ in this case $[k]=[K]$ is equivalent to
$$\calt_{M,\text{can}}(1)-\frac{\lambda(M)}{|H|}=\frac{K^2+s}{8}+DP.$$
This identity was verified in \cite{SWII}. In the sequel we verify
the corresponding identity valid for an arbitrary
$\spin^c$--structure.\end{h}

\begin{h}\textbf{The identity $\ssw^{OSz}_{M,[k]}=\ssw^{TCW}_{M,[k]}$}\qua
\label{sei15} Using
the above facts, this identity is equivalent  to
\begin{equation*}
\et_{M,[k]}(1)-\frac{\lambda(M)}{|H|}=\chi(HF^+(-M,[k]))-\min \chi_{k_r}+\frac{k_r^2+s}{8}.
\tag{$*$}\end{equation*}
We fix a $\spin^c$--structure  $[k]$ which corresponds to $\{a_0;a_1,\ldots,a_\nu\}$
as above (see Corollary~\ref{sei7}).
For the completeness of the description, we also need (see \cite{SWII}):
$$\hat{P}_1(t):=\frac{(t^\alpha-1)^{\nu-2}}{\prod_l (t^{\alpha/\alpha_l}-1)}.$$
Then by \cite{SWII} (see also Sections \ref{4.1} and \ref{l11}),
and with the notation  $\widetilde{a}=a_0+\sum_la_l/\alpha_l$ one has
$$\et_{M,[k]}(1)=\lim_{t\to 1} \Big(
\sum_{i\geq \widetilde{a}/e}\max\Big\{ 0,
1+a_0-ie_0+\sum_l\Big[\frac{-i\omega_l+a_l}{\alpha_l}\Big]
\Big\}t^{oi+\alpha\widetilde{a}}-\frac{1}{|H|}\hat{P}_1(t)\, \Big).$$
Using the inequalities (SI)  (see Proposition~\ref{sei5}) one deduce
that the contribution corresponding to the indices $i<0$ is zero.
Then using $\max\{0,x\}-\max\{0,-x\}=x$,
we get that
$$\et_{M,[k]}(1)-\chi(HF^+(-M,[k]))+\min \chi_{k_r}$$
equals the limit  $$ L:=\lim_{t\to 1} \Big(
P_{[k]}(t)-\frac{1}{|H|}\hat{P}_1(t)\, \Big),$$ where $$
P_{[k]}(t):=\sum_{i\geq 0}\Big(
1+a_0-ie_0+\sum_l\Big[\frac{-i\omega_l+a_l}{\alpha_l}\Big]
\Big)t^{oi+\alpha\widetilde{a}}.$$ In particular, the identity $(*)$
is equivalent to $\lambda(M)/|H|+(k_r^2+s)/8=L$. This identity was
verified in \cite{SWII} for the canonical $\spin^c$--structure (see
Section~\ref{4.1}). Therefore, it is enough to check a relative
identity (the difference between the two identities corresponding
to an arbitrary $[k]$ and the canonical $[K]$). Since
$k_r^2+s=K^2+s-8\chi(l'_{[k]})$, this relative identity is
\begin{equation*}-\chi(l'_{[k]})= \lim_{t\to 1} \Big(
P_{[k]}(t)-P_{[K]}(t)
\Big). \tag{${\rm r*}$}\end{equation*} Notice that once this identity is
proved, it also provides a new formula for $\chi(l'_{[k]})$. In
order to verify ${\rm(r*)}$ we consider the Laurent expansions of
$P_{[k]}(t)$ with respect to $t^\xo-1$. They have a pole of order
two, and the above limit shows  that  the part with negative
exponents is independent of $[k]$. For $[K]$ it was determined in
\cite{SWII}, hence one has:
$$P_{[k]}(t)=\frac{-e}{(t^\xo-1)^2}+\frac{-e-e\ev/2}{t^\xo-1}+\ \mbox{higher order terms.}$$
In the sequel we write $P(t)\sim Q(t)$ if $\lim_{t\to 1}(P(t)-Q(t))=0$.

In the verification of ${\rm(r*)}$, we  consider the same inductive
steps as in Section~\ref{sei8}.
First, we verify ${\rm(r*)}$ for $l'_{[k]}=-a_0g_0$.  Then (using $\xo=-e\alpha$ as well) one gets
$$P_{[k]}(t)-P_{[K]}(t)=P_{[K]}(t)(t^{\alpha a_0}-1)+
\sum_{i\geq 0}a_0t^{\xo i+\alpha a_0}$$
$$\sim -(t^{\alpha a_0}-1)\Big(
\frac{e}{(t^\xo-1)^2}+\frac{e+e\ev/2}{t^\xo-1} \Big)+
\frac{a_0t^{\alpha a_0}}{1-t^\xo}\sim
\frac{a_0^2}{2e}+\frac{a_0}{2}(1+\ev).$$
This is compatible with $\chi(l'_{[-a_0g_0]})$ given in Section~\ref{sei8}.

Next, consider $[k]$ and $[k^-_l]$ as in Section~\ref{sei8}. Notice that
$$\Big[\frac{-i \omega_l+a_l}{\alpha_l}\Big]-
\Big[\frac{-i \omega_l+a_l-1}{\alpha_l}\Big]$$
is 1 if $-i\omega_l+a_l$ is divisible by $\alpha_l$, otherwise is zero.
The divisibility is satisfied if $i=\{\omega_l'a_l/\alpha_l\}\alpha_l +j\alpha_l$ for some  $j\geq 0$.
Therefore,
$$P_{[k]}(t)-P_{[k^-_l]}(t)=P_{[k]}(t)(1-t^{-\alpha/\alpha_l})+
\sum_{i\omega_l\equiv a_l}t^{\xo i+\alpha
\widetilde{a}-\alpha/\alpha_l}$$
$$\sim \frac{t^{\xo\alpha_l \{\frac{\omega'_la_l}{\alpha_l}\} +\alpha\widetilde{a}-\frac{\alpha}{\alpha_l}}}
{1-t^{\xo\alpha_l}}+(t^{-\frac{\alpha}{\alpha_l}}-1)\Big(
\frac{e}{(t^\xo-1)^2}+\frac{e+e\ev/2}{t^\xo-1} \Big) $$ $$
\sim\frac{1}{2}
+\frac{\ev}{2\alpha_l}-\frac{1}{2e\alpha_l^2}+\frac{\widetilde{a}}{e\alpha_l}-
\Big\{\frac{\omega'_la_l}{\alpha_l}\Big\}.$$ This agrees with
$\chi(l'_{[k^-_l]})-\chi(l'_{[k]})$ from Section~\ref{sei8}, hence  ${\rm(r*)}$
and $(*)$ follow by induction.\end{h}

\begin{cors}\label{sei16} Again, since
$\sum_{[k]}\et_{M,[k]}(1)=0$, one gets
$$\lambda(M)=\sum_{[k]\in \Spin^c(M)}
  \chi(HF^+(M,[k]))-\frac{d(M,[k])}{2}.$$\end{cors}

\end{document}

%% file: gtoutput.tex
\def\ifplaintex{\expandafter\ifx\csname documentclass\endcsname\relax}

\ifplaintex 
\else
\fi

\expandafter\ifx\csname beginpicture\endcsname\relax
\expandafter\ifx\csname documentclass\endcsname\relax
\input pictex \else
\input prepictex \input pictex \input postpictex \fi\fi

\def\gt{{\mathsurround=0pt\it $\cal G\mskip-2mu$eometry \&\ 
$\cal T\!\!$opology}}        

\def\gtp{{\mathsurround=0pt\it $\cal G\mskip-2mu$eometry \&\ 
$\cal T\!\!$opology $\cal P\!$ublications}}  

\def\lognumber#1{\def\thelognumber{#1}}
\def\volumenumber#1{\def\thevolumenumber{#1}}
\def\papernumber#1{\def\thepapernumber{#1}}
\def\volumeyear#1{\def\thevolumeyear{#1}}

\def\pagenumbers#1#2{\def\startpage{#1}\def\finishpage{#2}}
\def\published#1{\def\publishdate{#1}}
\def\proposed#1{\def\theproposer{#1}}
\def\seconded#1{\def\theseconders{#1}}
\def\received#1{\def\receiveddate{#1}}

\def\accepted#1{\def\accepteddate{#1}}
\def\asciititle#1{\def\theasciititle{#1}}
\def\covertitle#1{\def\thecovertitle{#1}}
\def\coverauthors#1{\def\thecoverauthors{#1}}
\def\asciiauthors#1{\def\theasciiauthors{#1}}
\def\asciiaddress#1{\def\theasciiaddress{#1}}
\def\asciiemail#1{\def\theasciiemail{#1}}

\long\def\asciiabstract#1{\long\def\theasciiabstract{#1}}
\def\asciikeywords#1{\def\theasciikeywords{#1}}

\def\shorttitle#1{\def\theshorttitle{#1}}

\let\\\par\let\thelognumber\relax
\let\thevolumenumber\relax\let\thepapernumber\relax
\let\thevolumeyear\relax\let\thesamplenumber\relax\let\startpage\relax
\let\finishpage\relax\let\publishdate\relax\let\receiveddate\relax
\let\reviseddate\relax\let\accepteddate\relax\let\theasciititle\relax
\let\thecovertitle\relax\let\theasciiauthors\relax\let\theasciiaddress\relax
\let\theasciiabstract\relax\let\theasciikeywords\relax
\let\theasciiemail\relax\let\theshortauthors\relax\let\theshorttitle\relax
\let\thecoverauthors\relax

\long\def\maketitlep{   

\count0=\startpage

\gt\hfill      
\beginpicture
\setcoordinatesystem units <0.33truein, 0.33truein> point at 2.2 0.9
\setplotsymbol ({$\cal G$})
\plotsymbolspacing=9truept
\circulararc 315 degrees from 0 1 center at 0 0
\setplotsymbol ({$\cal T$})
\circulararc 315 degrees from 1 -1 center at 1 0
\endpicture
\break
{\small\ifx\thesamplenumber\relax 
Volume \else Sample
\fi\thevolumenumber\ (\thevolumeyear)
\startpage--\finishpage\nl
Published: \publishdate}
\vglue 0.5truein plus 0.4fil minus 0.1truein

{\parskip=0pt\leftskip 0pt plus 1fil\def\\{\par\smallskip}{\ifplaintex\large
\else\Large\fi\bf\thetitle}\par\medskip}   

\vglue 0pt plus 0.1fil 

{\parskip=0pt\leftskip 0pt plus 1fil\def\\{\par}{\sc\theauthors}
\par\medskip}

\vglue 0pt plus 0.1fil 

{\small\parskip=0pt\let\newline\\
{\leftskip 0pt plus 1fil\def\\{\par}{\sl\theaddress}\par}
\expandafter\ifx\theemail\relax    
\relax\else\vglue 5pt plus 0.02fil minus 2pt\def\\{\stdspace{\rm 
and}\stdspace} 
\cl{Email:\stdspace\tt\theemail}\fi
\ifx\theurl\relax                  
\relax\else\vglue 5pt plus 0.02fil minus 2pt\def\\{\stdspace{\rm 
and}\stdspace}
\cl{URL:\stdspace\tt\theurl}\fi\par}

\vglue 7pt plus 0.3fil minus 3pt

{\bf Abstract}
\vglue 5pt plus 0.1fil minus 2pt

\theabstract

\vglue 7pt plus 0.3fil minus 3pt

{\bf AMS Classification numbers}\quad Primary:\quad \theprimaryclass

Secondary:\quad \thesecondaryclass

\vglue 5pt plus 0.3fil minus 2pt

{\bf Keywords:}\quad \thekeywords

\vglue 10pt plus 0.5fil minus 5pt

{\small  Proposed: \theproposer\hfill Received: \receiveddate\nl
Seconded: \theseconders\hfill 
\ifx\reviseddate\relax                         
Accepted: \accepteddate                        
\else
Revised: \reviseddate                          
\fi}
\eject
}       

\let\maketitlepage\maketitlep
\let\maketitle\maketitlepage

\font\phead=cmsl9 scaled 950
\font=cmsl9 scaled 1050
\font\pnum=cmbx10 scaled 913
\font\lnum=cmbx10 
\font\pfoot=cmsl9 scaled 950
\font=cmsl9 scaled 1050
\ifplaintex
\headline{\vbox to 0pt{\vskip -4.5mm\line{\small\phead\ifnum
\count0=\startpage ISSN 1364-0380 (on line)
1465-3060 (printed) \hfill {\pnum\folio}\else\ifodd\count0\def\\{ }%
\ifx\theshorttitle\relax\thetitle\else\theshorttitle\fi\hfill{\pnum\folio}
\else\def\\{ and }{\pnum\folio}\hfill\ifx\theshortauthors\relax\theauthors
\else\theshortauthors\fi\fi\fi}\vss}}
\footline{\vbox to 0pt{\vglue 0mm\line{\small\pfoot\ifnum\count0=\startpage
\copyright\ \gtp\hfill\else
\gt, Volume \thevolumenumber\ (\thevolumeyear)\hfill\fi}\vss
}}
\else
\makeatletter
\def\@oddhead{{\small\lhead\ifnum\count0=\startpage ISSN 1364-0380 (on line)
1465-3060 (printed) \hfill {\lnum\number\count0}\else\ifodd\count0
\def\\{ }\ifx\theshorttitle\relax \thetitle \else\theshorttitle\fi\hfill
{\lnum\number\count0}\else\def\\{ and }{\lnum\number\count0}
\hfill\ifx\theshortauthors\relax 
\theauthors\else\theshortauthors\fi\fi\fi}}\def\@evenhead{\@oddhead}
\def\@oddfoot{\small\lfoot\ifnum\count0=\startpage\copyright\ \gtp\hfill\else
\gt, Volume \thevolumenumber\ (\thevolumeyear)\hfill\fi}
\def\@evenfoot{\@oddfoot}
\makeatother
\fi

\newwrite\gtoutfile
\long\gdef\makeheadfile{  
{\def\\{, }\def\s{ }
\immediate\openout\gtoutfile head.xxx
\immediate\write\gtoutfile{Proxy-for: \ifx\theasciiauthors\relax
\theauthors\else\theasciiauthors\fi\s<\ifx\theasciiemail\relax\theemail\else\theasciiemail\fi>}
\immediate\write\gtoutfile{\noexpand\\}
\immediate\write\gtoutfile{Authors: \ifx\theasciiauthors\relax
\theauthors\else\theasciiauthors\fi}
{\def\\{ }\immediate\write\gtoutfile{Title: \ifx\theasciititle\relax
\thetitle\else\theasciititle\fi}}
\immediate\write\gtoutfile{Subj-class: GT or SG or MG etc}
\immediate\write\gtoutfile{MSC-class: \theprimaryclass\ifx\thesecondaryclass\relax\else, \thesecondaryclass\fi}
\immediate\write\gtoutfile{Journal-ref: Geom. Topol. \thevolumenumber
(\thevolumeyear) \startpage-\finishpage}
\immediate\write\gtoutfile{Comments: Published by Geometry and Topology at}
\immediate\write\gtoutfile{\s\s http://www.maths.warwick.ac.uk/gt/GTVol\thevolumenumber/paper\thepapernumber.abs.html}
\immediate\write\gtoutfile{\noexpand\\}
\immediate\write\gtoutfile{}
\ifx\theasciiabstract\relax
\immediate\write\gtoutfile{\theabstract}\else
\immediate\write\gtoutfile{\theasciiabstract}\fi
\immediate\write\gtoutfile{}
\immediate\write\gtoutfile{\noexpand\\}
\immediate\write\gtoutfile{}
\immediate\closeout\gtoutfile}}  

\def\maketitlepage{\maketitlep\makeheadfile}
\let\maketitle\maketitlepage

%% file: 2005-23.bbl
\begin{thebibliography}{99}

\bibitem{Artin62} \textbf{M Artin},
\emph{Some numerical criteria for contractibility of curves on
algebraic surfaces}, Amer. J. Math. 84 (1962) 485--496
  \MR{0146182}

\bibitem{Artin66} \textbf{M Artin},
\emph{On isolated rational singularities of surfaces}, Amer. J. 
Math. 88 (1966) 129--136
  \MR{0199191}

\bibitem{GS} \textbf{R\,E Gompf}, \textbf{I\,A Stipsicz},
\emph{An Introduction to $4$--Manifolds and Kirby Calculus},
Graduate Studies in Mathematics 20, Amer. Math. Soc. (1999)
  \MR{1707327}

\bibitem{GRa} \textbf{H Grauert},
\emph{\"Uber Modifikationen  und exceptionelle
analytische Mengen}, Math. Ann. 146 (1962) 331--368
  \MR{0137127}

\bibitem{Laufer72} \textbf{H\,B Laufer}, \emph{On rational singularities},
Amer. J. Math. 94 (1972) 597--608
  \MR{0330500}

\bibitem{Laufer77} \textbf{H\,B Laufer},
\emph{On minimally elliptic singularities},
Amer. J. Math. 99 (1977) 1257--1295
  \MR{0568898}

\bibitem{Lescop} \textbf{C Lescop}, \emph{Global Surgery Formula for the
Casson--Walker
Invariant}, Annals of Mathematics Studies 140, Princeton
University Press (1996)
  \MR{1372947}

\bibitem{LMN} \textbf{I Luengo}, \textbf{A Melle-Hern\'andez},
\textbf{A N\'emethi},
\emph{Links and analytic invariants of superisolated
singularities}, J. Algebraic Geom. 2 (1993) 1--23
  \MR{1185605}

\bibitem{Ninv}  \textbf{A N\'emethi},
\emph{``Weakly'' Elliptic Gorenstein Singularities of Surfaces},
Invent. Math. 137 (1999) 145--167
  \MR{1703331}

\bibitem{Five} \textbf{A N\'emethi},
\emph{Five lectures on normal surface singularities}, from: ``Low
dimensional topology (Eger 1996/Budapest 1998)'', 
Bolyai Soc. Math. Stud. 8 (1999) 269--351
  \MR{1747271}

\bibitem{NL} \textbf{A N\'emethi},
\emph{Line bundles associated with normal surface singularities},
\arxiv{math.AG/0310084}

\bibitem{SWI} \textbf{A N\'emethi}, \textbf{L\,I Nicolaescu},
\emph{Seiberg--Witten invariants and surface singularities},
\gtref{6}{2002}{9}{269}{328}
  \MR{1914570}

\bibitem{SWII} \textbf{A N\'emethi}, \textbf{L\,I  Nicolaescu},
\emph{Seiberg--Witten invariants and surface singularities II
(singularities with good $\C^*$--action)},
J. London Math. Soc. 69 (2004) 593--607
  \MR{2050035}

\bibitem{NP} \textbf{W\,D Neumann},
\emph{A calculus for plumbing applied to the
topology of complex surface singularities and degenerating complex
curves},
Trans. Amer. Math. Soc. 268 (1981) 299--344
  \MR{0632532}

\bibitem{Neu} \textbf{W Neumann},
\emph{Abelian covers of quasihomogeneous surface
singularities}, from: ``Singularities Part 2 (Arcata Calif. 1981)'',
Proc. Sympos. Pure Math. 40 (1983) 233--244
  \MR{0713252}

\bibitem{Nico5} \textbf{L\,I Nicolaescu},
\emph{Seiberg--Witten invariants of rational homology 3--spheres},
Commun. Contemp. Math. 6 (2004) 833--866
  \MR{2111431}

\bibitem{OSz} \textbf{P\,S Ozsv\'ath}, \textbf{Z Szab\'o},
\emph{Holomorphic discs and topological invariants for closed three-spheres},
Ann. of Math. (2) 159 (2004) 1027--1158
\MR{2113019}

\bibitem{OSz7} \textbf{P\,S Ozsv\'ath}, \textbf{Z Szab\'o},
\emph{Holomorphic discs and three-manifold invariants: properties and
applications},
Ann. of Math. (2) 159 (2004) 1159--1245
  \MR{2113020}

\bibitem{OSzTr} \textbf{P\,S Ozsv\'ath}, \textbf{Z Szab\'o},
\emph{Holomorphic triangle
invariants and the topology of symplectic four-manifolds},
Duke Math. J. 121 (2004) 1--34
  \MR{2031164}

\bibitem{OSzAB} \textbf{P\,S Ozsv\'ath}, \textbf{Z Szab\'o},
\emph{Absolutely graded Floer homologies and intersection forms for
four-manifolds with boundaries},
Adv. Math. 173 (2003) 179--261
  \MR{1957829}

\bibitem{OSzP} \textbf{P\,S Ozsv\'ath}, \textbf{Z Szab\'o},
\emph{On the Floer homology of plumbed three-manifolds},
\gtref{7}{2003}{5}{185}{224}
  \MR{1988284}

\bibitem{Pinkh} \textbf{H Pinkham},
\emph{Normal surface singularities with ${\C}^*$ action},
Math. Ann. 117 (1977) 183--193
  \MR{0432636}

\bibitem{Stevens} \textbf{J Stevens},
\emph{Elliptic surface singularities and smoothing of curves},
Math. Ann. 267 (1984) 239--249
  \MR{0738250}

\bibitem{Tu5} \textbf{V\,G Turaev},
\emph{Torsion invariants of $\Spin^c$--structures on $3$--manifolds},
Math. Res. Lett. 4 (1997) 679--695
  \MR{1484699}

\end{thebibliography}
